\RequirePackage{amsthm}
\documentclass[sn-mathphys-num]{sn-jnl}

\usepackage{graphicx}%
\usepackage{multirow}%
\usepackage{amsmath,amssymb,amsfonts}%
\usepackage{lmodern}
\usepackage{mathabx}
\usepackage{mathrsfs}%
\usepackage[title]{appendix}%
\usepackage{xcolor}%
\usepackage{textcomp}%
\usepackage{manyfoot}%
\usepackage{booktabs}%
\usepackage{algorithm}%
\usepackage{algorithmicx}%
\usepackage{algpseudocode}%
\usepackage{listings}%
\usepackage{enumitem}
\usepackage{graphicx}
\usepackage{subfig}
\usepackage[misc]{ifsym}
\usepackage{lineno}

    \makeatletter
\def\@fnsymbol#1{\ensuremath{\ifcase#1\or
   \mathsection\or \mathparagraph\or \|\or **\or \dagger\dagger
   \or \ddagger\ddagger \else\@ctrerr\fi}}
    \makeatother

\newcommand{\lv}[0]{\mathbf{l}}
\newcommand{\uv}[0]{\mathbf{u}}
\newcommand{\vv}[0]{\mathbf{v}}
\newcommand{\qv}[0]{\mathbf{q}}
\newcommand{\xv}[0]{\mathbf{x}}
\newcommand{\yv}[0]{\mathbf{y}}
\newcommand{\pv}[0]{\mathbf{p}}

\newcommand{\iprod}[2]{\left\langle {#1}, {#2} \right\rangle}
\newcommand{\norm}[1]{\left\lVert#1\right\rVert}

\newcommand{\dv}[0]{\mathbf{d}}
\newcommand{\ev}[0]{\mathbf{e}}

\newcommand{\cv}[0]{\mathbf{c}}
\newcommand{\bv}[0]{\mathbf{b}}
\newcommand{\hv}[0]{\mathbf{h}}
\newcommand{\gv}[0]{\mathbf{g}}
\newcommand{\sv}[0]{\mathbf{s}}
\newcommand{\wv}[0]{\mathbf{w}}
\newcommand{\rv}[0]{\mathbf{r}}

\newcommand{\kv}[0]{\mathbf{k}}
\newcommand{\tv}[0]{\mathbf{t}}
\newcommand{\zv}[0]{\mathbf{z}}

\newcommand{\conv}[0]{{\mathbf{conv}}}
\newcommand{\convc}[0]{{\overline{\mathbf{conv}}}}

\DeclareMathOperator{\epi}{\operatorname{epi}}

\DeclareMathOperator{\Rd}{\operatorname{R}}
\DeclareMathOperator{\OPT}{\operatorname{OPT}}
\DeclareMathOperator{\DUAL}{\operatorname{DUAL}}

\newcommand{\E}{\mathbb{E}}

\newcommand{\R}{\mathbb{R}}

\newcommand{\blue}[1]{\textcolor{black}{#1}}

\newcommand{\interior}{\textbf{Int}}

\newcommand*\xbar[1]{%
  \hbox{%
    \vbox{%
      \hrule height 0.5pt 
      \kern0.5ex
      \hbox{%
        \kern-0.1em
        \ensuremath{#1}%
        \kern-0.1em
      }%
    }%
  }%
}

\theoremstyle{thmstyleone}%
\newtheorem{theorem}{Theorem}
\theoremstyle{thmstyletwo}%
\theoremstyle{thmstylethree}%
\newtheorem{example}{Example}%
\newtheorem{remark}{Remark}%
\newtheorem{corollary}{Corollary}%

\newtheorem{definition}{Definition}%
\newtheorem{asp}{Assumption}%
\newtheorem{lemma}{Lemma}
\newtheorem{claim}{Claim}

\begin{document}

\title{Asymptotically Strong Dual of Smooth Nonconvex Problems: A Smooth Shapley–Folkman Perspective}




\author[1]{\fnm{Santanu} \sur{S. Dey}}\email{sdey30@gatech.edu}

\author*[1]{\fnm{Jingye} \sur{Xu}}\email{jxu673@gatech.edu}

\affil[1]{\orgdiv{H. Milton Stewart School of Industrial and Systems Engineering}, \orgname{Georgia Institute of Technology}, \orgaddress{\street{ 755 Ferst Dr NW}, \city{Atlanta}, \postcode{30332}, \state{GA}, \country{USA}}}




\abstract{The Shapley–Folkman Lemma, a foundational result in convex geometry,  asserts that the nonconvexity of a Minkowski sum of $n$-dimensional bounded nonconvex sets does not accumulate once the number of summands exceeds the dimension $n$, and thus the sum becomes approximately convex. 
This Lemma, originally published by Starr is used to show quasi-equilibrium in nonconvex market models in economics and equivalently used to estimate the duality gap of the Lagrangian dual of separable nonconvex problems.
Given its foundational nature, we pose the following geometric question: Is it possible for the nonconvexity of the Minkowski sum of 
$n$-dimensional nonconvex sets to vanish asymptotically, rather than merely converge to a constant, as the number of summands increases?
We answer this question affirmatively. 
First, we provide an elementary geometric proof of the Shapley–Folkman Lemma based on the facial structure of the convex hull of each set under a pointedness assumption. 
This leads to a small refinement on the classical error bound derived from the Lemma in the pointed case.
Building on this new geometric perspective, we next show that if the sets satisfy a certain local smoothness condition, which naturally arises, for example, in the epigraphs of smooth functions, then their Minkowski sums converge directly to a convex set as the nonconvexity measure vanishes asymptotically.
This result has direct implications for optimization and economic theory. In particular, we show that the Lagrangian dual of block-structured smooth nonconvex problems, possibly subject to additional sparsity constraints, is asymptotically tight under mild assumptions. 
This sharply improves upon classical Shapley–Folkman based bounds, which generally predict a nonvanishing duality gap. 
Equivalently, from an economic perspective, this result implies that, asymptotically as the size of the economy increases, competitive equilibrium exists with smooth utility functions, notwithstanding underlying nonconvexities.
}




\maketitle

\section{Introduction}

Consider a block-structured nonconvex program:
\begin{subequations}
\label{prob}
\begin{align}
\OPT(\bv) := & \inf \sum_{i=1}^k \iprod{\cv^{(i)}}{\xv^{(i)}} \\
        & \text{s.t. } \sum_{i=1}^k B^{(i)} \xv^{(i)} \leq \bv, \label{eq:coupling} \\
        & \ \ \ \ \ \ \xv^{(i)} \in \mathcal{X}^{(i)} \subseteq \R^n,\forall i \in [k],
\end{align}
\end{subequations}
where $\mathcal{X}^{(i)}$ is a closed possibly nonconvex set of $\R^n$,  $B^{(i)} \in \R^{m \times n}$ and $\bv \in \R^{m}$. 
Although the objective and the coupling constraints in (\ref{prob}) are linear, this formulation is general enough to capture many nonlinearities, by embedding these into the feasible sets $\mathcal{X}^{(i)}$, potentially through the introduction of suitable auxiliary variables.

An important relaxation of (\ref{prob}) is the Lagrangian dual 
of (\ref{prob}), obtained by dualizing the coupling constraints (\ref{eq:coupling}) across various blocks:
\begin{equation}
    \begin{aligned}
        L(\lambda) := & \inf_{\xv} \left(\sum_{i=1}^k \iprod{\cv^{(i)}}{\xv^{(i)}} + \iprod{\lambda}{B^{(i)} \xv^{(i)}} \right) - \iprod{\lambda}{\bv} \\
        & \text{s.t. } \xv^{(i)} \in \mathcal{X}^{(i)} \subseteq \R^n,\forall i \in [k],
    \end{aligned}
\end{equation}
and
\begin{equation}
\begin{aligned}
\label{DUAL}
\DUAL(\bv) := \sup_{\lambda \geq 0} L(\lambda).
\end{aligned}
\end{equation}
The Lagrangian dual of the block-structured problem~(1) plays a central role in
both optimization and economic theory.
From an optimization perspective, dualizing the coupling constraints yields a
convex (possibly nonsmooth) program (\ref{DUAL}) regardless of the nonconvexity of the feasible sets $\mathcal{X}^{(i)}$, thereby providing a tractable and globally
valid lower bound on $\OPT(\mathbf{b})$.
Moreover, the dual formulation exploits the separable block structure of the
primal problem: for any fixed multiplier $\lambda \ge 0$, the Lagrangian
$L(\lambda)$ decomposes into independent subproblems across blocks, enabling
scalable algorithms based on decomposition, subgradient methods, and
price-based coordination \cite{lubbecke2010column,desaulniers2006column,cifuentes2025lagrangian}. 
From an economic viewpoint, the dual variables $\lambda$ admit a natural
interpretation as equilibrium prices associated with the shared resource
constraint $\sum_{i=1}^k B^{(i)} x^{(i)} \le \mathbf{b}$ \cite{starr1969quasi}.
In this interpretation, each block represents an individual agent or firm that
optimizes its own objective given prices $\lambda$, while the dual objective
aggregates these decentralized decisions and enforces market clearing.
Consequently, the duality gap between $\OPT(\mathbf{b})$ and
$\DUAL(\mathbf{b})$ quantifies the extent to which a competitive equilibrium
exists or fails to exist in the presence of nonconvexities, making the
Lagrangian framework a unifying bridge between large-scale optimization and
equilibrium analysis in economic theory.

One of the most desirable properties one aims to establish for the Lagrangian dual, in both optimization and economic theory, is \emph{strong duality}.
In optimization, strong duality certifies that the optimal value of a (possibly relaxed) dual problem coincides with that of the primal, thereby justifying convex relaxations, decomposition methods, and price-based algorithms. In economic theory, strong duality corresponds to the existence of supporting prices that decentralize a socially optimal allocation, that is, establishing the existence of a competitive equilibrium.

Many textbooks state that strong duality, that is $\OPT(\bv) = \DUAL(\bv)$, holds under certain regularity conditions when all $\mathcal{X}^{(i)}$s are convex.
However, strong duality actually  only requires the convexity of \blue{the} linear image of (\ref{prob}) onto the space of objective and constraints~\cite{ben2001lectures}. More specifically, let

\begin{align*}
    & \mathcal{P}^{(i)} := \left\{ (t,\dv) \;\middle\vert\;
   \begin{array}{@{}l@{}} \exists \xv \in \mathcal{X}^{(i)}, t= \iprod{\cv^{(i)}}{\xv^{(i)}},\dv = B^{(i)} \xv^{(i)}
   \end{array}
   \right\}, \\
   & \mathcal{P} := \sum_{i=1}^k \mathcal{P}^{(i)},
\end{align*}
where we overload standard sum for vectors with Minkowski sum for sets and $\mathcal{P}$ is the Minkowski sum of all $\mathcal{P}^{(i)}$s.
Under suitable regularity conditions, {strong duality} can be guaranteed whenever the image set $\mathcal{P}$ is convex~\cite{ben2001lectures}.
In particular, when each set $\mathcal{X}^{(i)}$ is convex, the overall set $\mathcal{P}$ remains convex because both linear mappings and Minkowski sums preserve convexity. 
Therefore, the condition frequently shown in textbooks can be viewed as a sufficient condition.
Although this distinction may appear subtle, it is in fact crucial---it underpins several nontrivial results concerning the {strong duality} in Lagrangian relaxations of nonconvex programs~\cite{polik2007survey,chandrasekaran2025lagrangian}. 
A notable example is the celebrated $S$-Lemma~\cite{polik2007survey}, which establishes that the Lagrangian dual of a quadratic program with a single quadratic constraint is tight. 
This result can be interpreted as a direct consequence of Dines Lemma \cite{dines1941mapping} that the image of two homogeneous quadratic mappings is always convex.

In this paper, we adopt this perspective and investigate the convexity of $\mathcal{P}$ to ensure the strong duality of the Lagrangian dual of~\eqref{prob}. 
While each individual set $\mathcal{P}^{(i)}$ may be nonconvex, it is possible that their Minkowski sum $\mathcal{P}$ closely approximates its convex hull. 
This phenomenon is formally characterized by the \emph{Shapley--Folkman Lemma}~\cite{starr1969quasi}. 
In the setting where $k \ge \blue{m+1}$, the Shapley-Folkman Lemma states that
\begin{align*}
    \forall \xv \in \conv(\mathcal{P}), \; \exists I \subseteq [k], |I| \leq \blue{m+1} 
    \text{ such that } 
    \xv \in 
    \left( \sum_{i \in [k] \setminus I} \mathcal{P}^{(i)} \right)
    + 
    \left( \sum_{i \in I} \conv\left(\mathcal{P}^{(i)}\right) \right).
\end{align*}

Here $[k] = \{1, \dots, k\}$. 
Intuitively, the Shapley-Folkman Lemma implies that every point in $\conv(\mathcal{P})$ is close to some point in $\mathcal{P}$, up to at most \blue{$m+1$} `nonconvex summands'. 
If all $\mathcal{P}^{(i)}$s are contained within a ball of radius $R$, then a quantitative measure of this approximate convexity can be obtained. 
In particular, under the Hausdorff distance with $l_2$ norm, it has been shown~\cite{fradelizi2018convexification} that
\begin{align*}
    d_H(\mathcal{P}, \conv(\mathcal{P})) \le \sqrt{\blue{m+1}}R.
\end{align*}
This quantitative bound can be employed to derive an {a priori} estimate of the duality gap 
$\Delta := \OPT(\bv) - \DUAL(\bv)$ associated with problem~\eqref{prob} (see Theorem \ref{thm_duality_gap} in Section~\ref{sec:SF2gap}). 


The Shapley–Folkman lemma was originally introduced by Starr~\cite{starr1969quasi} in the study of quasi-equilibria in nonconvex market models in economics. In this context, it can be interpreted as an approximate form of strong duality for Lagrangian duals. Its connection to optimization was later developed in works such as \cite{bertsekas1983optimal, aubin1976estimates, udell2016bounding}, which used the lemma to quantify the duality gap of problem \eqref{prob}.
More recently, several works have focused on refining the quantitative bounds provided by the Shapley–Folkman lemma, particularly on improving estimates of the Hausdorff distance 
$d_H(\mathcal{P}, \conv(\mathcal{P}))$ or similar notion of nonconvexity; see, for example, \cite{bi2016refined, kerdreux2023stable, nguyen2022near}. An analog of the Shapley–Folkman lemma for discrete convex sets has also been developed in \cite{murota2025shapley}.


In this work, we seek to refine the quantitative characterization of the Hausdorff distance 
$d_H(\mathcal{P}, \conv(\mathcal{P}))$, given its central role in assessing the tightness 
of Lagrangian duals. We begin by presenting a new and elementary geometric proof of the 
Shapley-Folkman Lemma, which leverages the facial structure of each 
$\operatorname{conv}(\mathcal{P}^{(i)})$. This perspective enables us to refine the classical 
error bound under a pointedness assumption (see Theorem~\ref{thm:generalrefinement} in Section~\ref{sec:newbnd}). 
Building on this geometric insight, we further demonstrate that the Minkowski sum of the sets 
$\mathcal{P}^{(i)}$ can directly converge to its convex hull as the number of summands grows, 
i.e., $d_H(\mathcal{P}, \operatorname{conv}(\mathcal{P})) \to 0$ as $k \to \infty$, 
provided enough sets satisfy some mild ``local smoothness'' condition (see Theorem \ref{thm_rem_convex_ball} in Section~\ref{sec:newbnd}). To the best of our knowledge, this is the first result establishing an asymptotically vanishing nonconvexity measure for Minkowski sums of nonconvex sets within the framework of the Shapley–Folkman Lemma.

Such smoothness naturally arises 
in the epigraphs of smooth (possibly nonconvex) functions, allowing us to establish that the following problem 
admits an asymptotically tight Lagrangian dual as $k \to \infty$:
\begin{equation}
    \label{prob_smooth_sparsity}
    \begin{aligned}
         \inf_{\xv} \quad & \sum_{i=1}^{k} f^{(i)}(\xv^{(i)}) \\
    \text{s.t.} \quad & \sum_{i=1}^k B^{(i)} \xv^{(i)} \leq \bv, \\
    & \|\xv^{(i)}\|_0 \le s^{(i)}, \quad \forall i \in [k],
    \end{aligned}
\end{equation}
where each $f^{(i)}$ is a smooth (possibly nonconvex) function, $B^{(i)}$s satisfy some mild conditions, and $s^{(i)}$ is not too small, i.e., $s^{(i)} \geq m,\forall i \in [k]$ where  $B^{(i)}$ has $m$ rows (see {Theorem} \ref{thm_application} in Section~\ref{sec:app}) . To model problems without sparsity constraints, one may simply drop these constraints, that is let $s^{(i)} = \textup{number of variables in }\xv^{i}$, for all $i \in [k]$. 

Finally, we note here that this work focuses on Minkowski sum of $\mathcal{P}^{(i)}$s without averaging since its convexity directly reflects the duality gap of problem~(\ref{prob}). In the classical applications 
of the Shapley--Folkman Lemma, the analysis is typically performed on the Minkowski {average} 
of the sets $\mathcal{P}^{(i)}$. 
This choice arises because the error bound 
$d_H(\mathcal{P}, \conv(\mathcal{P}))$ does not vanish in the classic analysis of Shapley-Folkman Lemma but $d_H(\mathcal{P}, \conv(\mathcal{P}))$ is both homogeneous and independent of $k$, 
which immediately yields 
\[
d_H\!\left(\tfrac{1}{k}\mathcal{P}, \tfrac{1}{k}\conv(\mathcal{P})\right) 
\le \tfrac{\sqrt{\blue{m+1}}R}{k},
\]
implying a linear convergence of the averaged sum toward its convex hull as $k$ increases. 
In contrast, our analysis applies directly to the unnormalized Minkowski sum. When averaging is applied, our result demonstrates 
a {superlinear} convergence rate, thereby strengthening the classical result.

The remainder of the paper is organized as follows. In Section~\ref{sec:pre}, we present standard notation and preliminary results. Section~\ref{sec:SF2gap} establishes a connection between $d_H(\mathcal{P}, \conv(\mathcal{P}))$ and $ \Delta := \OPT(\bv) - \DUAL(\bv)$ in Theorem~\ref{thm_duality_gap}.  Section~\ref{sec:newbnd} presents a new proof of Shapley-Folkman Lemma under the assumption of pointedness of the convex hull of the non-convex sets (see Lemma~\ref{lem:face_decompose} and discussion following it) which allows us to establish a slightly better bound on $d_H(\mathcal{P}, \conv(\mathcal{P}))$, in  Theorem~\ref{thm:generalrefinement}. Our main result of showing that $d_H(\mathcal{P}, \conv(\mathcal{P}))$ goes to zero as $k$ goes to $\infty$ under the presence of a local smoothness property, is presented in this section in Theorem~\ref{thm_rem_convex_ball}. Finally, in Section~\ref{sec:app}, we establish the strong duality result stated informally above, formalized as Theorem~\ref{thm_application}. 

\section{Preliminary Results and Notation}\label{sec:pre}

In this section, we introduce several useful notation and results to study the Shapley-Folkman Lemma and the Lagrangian dual of (\ref{prob}).

\subsection{Notation}

For a set \( A \subset \mathbb{R}^d \), we define its \emph{radius} as $
\Rd(A) := \inf_{x \in \mathbb{R}^d} \left\{ r \ge 0 : A \subseteq \mathcal{B}(x,r) \right\}$,
where \( \mathcal{B}(x,r) := \{ y \in \mathbb{R}^d : \|y - x\|_2 \le r \} \) denotes the closed Euclidean ball of radius \( r \) centered at \( x \).
Given two sets \( A, B \subset \mathbb{R}^d \), their \emph{Minkowski sum} is defined as
$
A + B := \{ a + b : a \in A,\; b \in B \}.
$
For any scalar \( \lambda \ge 0 \), the \emph{Minkowski scaling} of \( A \) is
$
\lambda A := \{ \lambda a : a \in A \}.
$We denote by
$
\conv(A)
:= \left\{ \sum_{i=1}^m \lambda_i a_i \;\middle|\;
a_i \in A,\ \lambda_i \ge 0,\ \sum_{i=1}^m \lambda_i = 1 \right\}
$
the convex hull of \( A \). The \emph{closed convex hull} is defined as
$
\overline{\conv}(A) := \operatorname{cl}(\conv(A)),
$
where \( \operatorname{cl}(A) \) denotes the topological closure of \( A \). The relative boundary and the recession cone of a closed convex set $C$ are denoted using $\textup{relbd}(C)$ and $\text{rec}(C)$. The kernel of a linear map $\mathcal{A}$ is denoted by $\ker(\mathcal{A})$. We use $\epi(f)$ to denote the epigraph of a function $f$. A set is called \emph{pointed} if it does not contain any affine line.

\subsection{Nonconvexity measure and convex hull operator}

For a given set $A \subseteq \R^n$, there exists a rich body of literature on how to measure non-convexity of $A$ due to its wide-ranging applications in nonconvex optimization. Throughout this paper, we will focus on two measures:
\begin{enumerate}
    \item Hausdorff distance from the convex hull: Given a closed convex set $K$ with $\textbf{0}$ in its interior, we define
    \begin{align*}
        \Phi^{K}(A) := \inf\left\{r \geq 0 : \conv(A) \subseteq A + r K  \right\}.
    \end{align*}
    When $K$ is omitted, we refer it as the Hausdorff distance from the convex hull under $l_2$ norm:
    \begin{align*}
        \Phi(A) & := \inf\left\{r \geq 0 : \conv(A) \subseteq A + r \mathcal{B}(\textbf{0},1)  \right\} \\
       & \blue{= \sup_{\xv \in \conv(A)} \inf_{\yv \in A} \norm{\xv-\yv}_2.}
    \end{align*}
    \item The inner radius of a nonconvex set~\cite{starr1969quasi}: 
    \begin{align*}
    \Xi(A) := \sup_{\xv \in \conv(A)} \left\{\inf{\Rd(T):T \subseteq A, \xv \in \conv(T})\right\}.
\end{align*}
\begin{figure}[!ht]
\centering
\begin{minipage}{.5\textwidth}
  \centering
  \includegraphics[width=.5\linewidth]{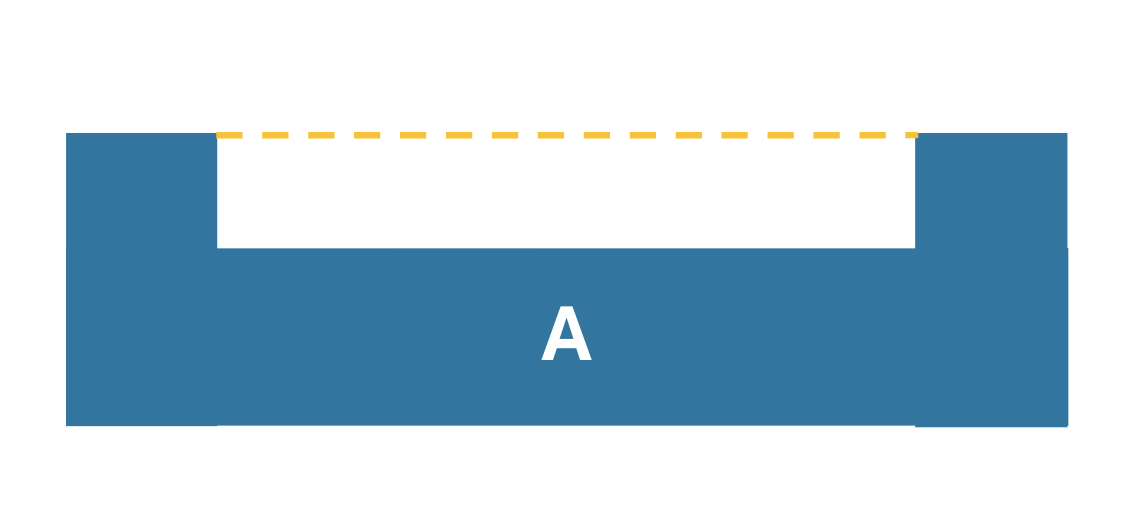}
  \captionof{figure}{$A$ and $\conv(A)$}
  \label{fig:AandConv}
\end{minipage}%
\begin{minipage}{.5\textwidth}
  \centering
  \includegraphics[width=.5\linewidth]{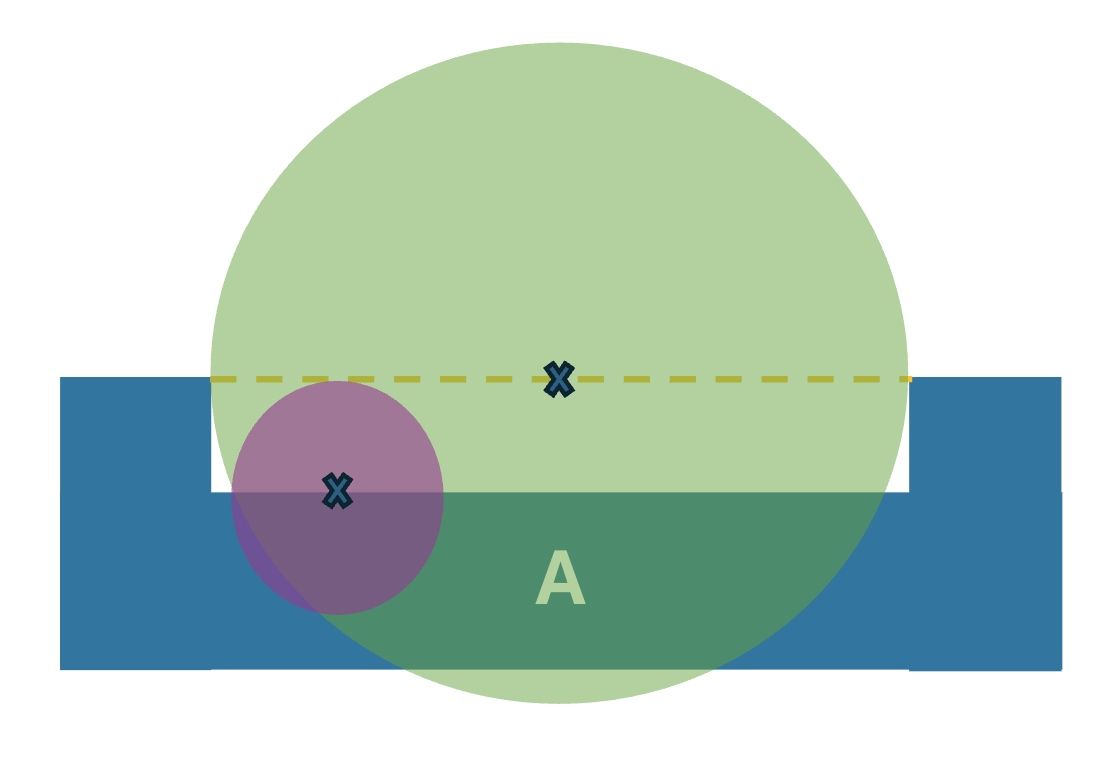}
  \captionof{figure}{two nonconvexity measures where \blue{the} pink ball corresponds to $\Phi(A)$ and \blue{the} green ball corresponds to $\Xi(A)$.}
  \label{fig:two_nonconvex_measure}
\end{minipage}
\end{figure}
\end{enumerate}


The measure $\Xi(\cdot)$, introduced in the seminal work by Starr \cite{starr1969quasi}, is used to establish the notion of quasi-equilibrium in nonconvex market models in economics. In this paper, we focus primarily on $\Phi(\cdot)$, as it is more directly related to the duality gap in Lagrangian dual (see Theorem \ref{thm_duality_gap}). See 
Figure\blue{s} \ref{fig:AandConv} and \ref{fig:two_nonconvex_measure} for illustration of the $\Phi(\cdot)$ and $\Xi(\cdot)$ nonconvexity measure.
Other nonconvexity measures, exhibiting distinct behaviors and intended for different analytical purposes, are not listed in this paper; for a comprehensive overview, we refer readers to the survey \cite{fradelizi2018convexification}. We present next a 
standard result on these non-convexity measures. 

\begin{lemma}
    \label{lem:relation_measure}
\cite{wegmann1980einige,fradelizi2018convexification} Let $A$ be a subset of $\R^n$, then
    \begin{align*}
        \Phi(A) \leq \Xi(A) \leq \Rd(A).
    \end{align*}
\end{lemma}




The analysis in this paper relies heavily on the convex hull operator $\conv(\cdot)$. Two structural properties of convex hulls are used repeatedly: commutativity with Minkowski sums and commutativity with linear projections. 
We restate these well-known facts below for convenience.



\begin{lemma}
\cite{fradelizi2018convexification}
\label{claim_convex_and_sum}
($\conv$ commutes with Minkowski sum)
Let $A_1,\dots,A_{k}$ be nonempty sets of $\R^n$. Then it follows that
\begin{align*}
    \conv\left(\sum_{i=1}^k A_i\right) = \sum_{i=1}^k \conv\left(A_i\right)
\end{align*}
\end{lemma}


\begin{lemma}
\cite{hiriart1993convex}
\label{claim_convex_and_proj}
($\conv$ commutes with linear projection)
\begin{align*}
   \conv\left(\mathcal{P}^{(i)}\right) & \blue{= \conv\left(\left\{ (t,\dv) \;\middle\vert\;
   \begin{array}{@{}l@{}} \exists \xv^{(i)} \in \mathcal{X}^{(i)}, t= \iprod{\cv^{(i)}}{\xv^{(i)}},\dv = B^{(i)} \xv^{(i)}
   \end{array}
   \right\}\right)}  \\ & = \left\{ (t,\dv) \;\middle\vert\;
   \begin{array}{@{}l@{}} \exists \xv^{(i)} \in \conv\left(\mathcal{X}^{(i)}\right), t= \iprod{\cv^{(i)}}{\xv^{(i)}},\dv = B^{(i)} \xv^{(i)}
   \end{array}
   \right\}.
\end{align*}
\end{lemma}



\subsection{Shapley-Folkman Lemma}

\begin{theorem}
    \label{thm_shapley_lemma}
    (Shapley-Folkman Lemma) Let $A_1,\dots,A_{k}$ be nonempty closed sets of $\R^n$ with $k \geq n+1$. Let $\xv \in \conv(\sum_{i\in[k]} A_{i} ) = \sum_{i \in [k]} \conv(A_i)$. Then there exists $\mathcal{I} \subseteq [k]$ with cardinality at most $n$ such that
    \begin{align*}
        \xv \in \sum_{i \in \mathcal{I}} \conv(A_i) + \sum_{i \in [k] \setminus \mathcal{I} } A_i.
    \end{align*}
\end{theorem}
\blue{The} Shapley--Folkman Lemma is a direct consequence of Carathéodory's theorem, and there are many different ways to prove it \cite{bi2016refined,kerdreux2023stable}. We defer our geometric proof to Lemma \ref{lem:face_decompose}, which is presented in Section~\ref{sec:newbnd}.

\blue{The} Shapley--Folkman Lemma can be interpreted as saying that the nonconvexity of the Minkowski sum of nonconvex sets does not continue to accumulate with the number of summands once it reaches \blue{a} certain threshold depending only on the dimension. This phenomenon can be quantified via $\Phi(\cdot)$ and $\Xi(\cdot)$ by the Corollary \ref{cor_quan_SFL}. In general, estimating the nonconvexity measures $\Phi(\cdot)$ and $\Xi(\cdot)$ is challenging. Therefore, many applications of the Shapley--Folkman Lemma consider bounded nonconvex sets and use the second part of Corollary \ref{cor_quan_SFL}.

\begin{corollary}
\label{cor_quan_SFL}
    Let $A_1,\dots,A_{k}$ be nonempty sets of $\R^n$ with $k \geq n$ such that $\Xi(A_i) \leq \beta,\forall i \in [k]$. Then it follows that
\begin{align*}
& \Phi\left(\sum_{i\in[k]} A_i\right) \leq \Xi\left(\sum_{i\in[k]} A_i\right)  \leq \sqrt{n} \beta  \text{ and therefore }
 \Phi\left(\frac{1}{k}\sum_{i\in[k]} A_i\right) \leq \frac{\sqrt{n}}{k} \beta.
 \end{align*}
 By Lemma \ref{lem:relation_measure}, if $\Rd(A_i) \leq \gamma,\forall i \in [k]$, then it follows that
 \begin{align*}
& \Phi\left(\sum_{i\in[k]} A_i\right) \leq \Xi\left(\sum_{i\in[k]} A_i\right)  \leq \sqrt{n} \gamma  \text{ and therefore }
 \Phi\left(\frac{1}{k}\sum_{i\in[k]} A_i\right) \leq \frac{\sqrt{n}}{k} \gamma.
 \end{align*}
\end{corollary}

\section{Duality gap of Lagrangian relaxation}\label{sec:SF2gap}

A key quantity associated with the Lagrangian dual that reflects the tightness of the relaxation is the duality gap (also called the additive gap) which is defined as:
\begin{align*}
    \Delta := \OPT(\bv) - \DUAL(\bv).
\end{align*}
In this section, we present a result showing how to upper bound $\Delta$ based on upper bounds on $\Phi(\mathcal{P})$.

To achieve this, we first use the following primal characterization~\cite{boland2018combining,geoffrion2009lagrangean} of Lagrangian dual (bi-conjugacy in convex analysis \cite{hiriart2004fundamentals}) to reduce $\DUAL(\bv)$ from min-max optimization to a minimization problem in the original space. Consider
\begin{equation}
    \label{eq:primal_char}
    \begin{aligned}
        \OPT_{\textbf{L}}(\bv) := & \inf \sum_{i=1}^k \iprod{\cv^{(i)}}{\xv^{(i)}} \\
        & \text{s.t. } \sum_{i=1}^k B^{(i)} \xv^{(i)} \leq \bv, \\
        & \ \ \ \ \ \ \xv^{(i)} \in \convc(\mathcal{X}^{(i)}) \subseteq \R^n,\forall i \in [k].
    \end{aligned}
\end{equation}


Under \blue{a} proper regularity condition \cite{geoffrion2009lagrangean,hiriart2004fundamentals}, it follows that $\OPT_{\textbf{L}}(\bv) = \DUAL(\bv)$.
Common regularity conditions include the existence of a Slater point or the assumption that each $\conv(\mathcal{X}^{(i)})$ is a polyhedron \cite{dey2025geoffrion,lemarechal2001geometric}. The absence of such regularity may result in cases where $\OPT_{\textbf{L}}(\bv) > \DUAL(\bv)$ \cite{lemarechal2001geometric}. Throughout this paper, we exclude these pathological instances and directly assume that the following equivalence holds:
\begin{asp}
    \label{asp:strong_dual}
    $\OPT_{\textbf{L}}(\bv) = \DUAL(\bv)$.
\end{asp}







%
We finally present the main result of this Section showing that the tightness of \blue{ the } Lagrangian dual of (\ref{prob}) directly relies on the convexity of $\mathcal{P} = \sum_{i=1}^k \mathcal{P}^{(i)}$ 
by quantifying $\Delta$ in terms of $\Phi(\mathcal{P})$. This theorem is similar to the gap formula in \cite{dubois2025frank}.


\begin{theorem}
    \label{thm_duality_gap}
    Let $\mathcal{E} :=\Phi\left(\mathcal{P}\right)$. Under Assumption \ref{asp:strong_dual} and $\conv(\mathcal{X}^{(i)})$ being closed for all $i \in [k]$, it follows that
    \begin{align*}
        & \OPT(\bv+\mathcal{E} \textbf{1}) - \mathcal{E} \leq \DUAL(\bv) \leq \OPT(\bv),
    \end{align*}
    and therefore
    \begin{align*}
                & \Delta \leq \mathcal{E} + \OPT(\bv) - \OPT(\bv+\mathcal{E} \textbf{1}).
    \end{align*}
    \begin{proof}
Under Assumption \ref{asp:strong_dual}, we know that $\DUAL(\bv) = \OPT_{\textbf{L}}(\bv)$. For an \blue{$\epsilon > 0$}, 
let \blue{$\xv_\textbf{L} = (\xv_\textbf{L}^{(1)},\cdots,\xv_\textbf{L}^{(k)})$} be \blue{an} $\epsilon$-optimal solution of (\ref{eq:primal_char}). By weak duality, we know that
\begin{equation}
    \label{eq:wd_1}
    \begin{aligned}
        \OPT_{\textbf{L}}(\bv) + \epsilon = \sum_{i=1}^k \iprod{\cv^{(i)}} {\xv^{(i)}_{\textbf{L}}} =  \DUAL(\bv) + \epsilon \leq \OPT(\bv) + \epsilon.
    \end{aligned}
\end{equation}
Let $\pv^{(i)}_\textbf{L} := \begin{bmatrix}
    \iprod{\cv^{(i)}} {\xv^{(i)}_{\textbf{L}}} \\
    B^{(i)} \xv^{(i)}_{\textbf{L}}
\end{bmatrix} \in \conv\left(\mathcal{P}^{(i)}\right)$ where the last containment is due to the fact that $\conv(\mathcal{X}^{(i)})$ is closed (which implies $\xv^{(i)}_{\textbf{L}} \in \conv(\mathcal{X}^{(i)})$) and using Lemma~\ref{claim_convex_and_proj}.
By Lemma \ref{claim_convex_and_sum}, \blue{$\sum_{i=1}^k \pv^{(i)}_{\textbf{L}} \in \sum_{i=1}^k \conv\left(\mathcal{P}^{(i)}\right) = \conv(\mathcal{P})$}. By the definition of $\Phi\left(\mathcal{P}\right)$, there exists some $\xv_{*} \in \mathcal{X}^{(i)}$ and $\pv^{(i)}_*:= \begin{bmatrix}
    \iprod{\cv^{(i)}} {\xv^{(i)}_*} \\
    B^{(i)} \xv^{(i)}_*
\end{bmatrix}$ such that
\begin{align}\label{eq:E}
    \norm{ \left( \sum_{i=1}^k \pv^{(i)}_\textbf{L} \right) - \left( \sum_{i=1}^k \pv^{(i)}_\textbf{*} \right) }_2 \leq \mathcal{E}.
\end{align}
This means that \blue{$\xv_{*}$} is \blue{a} feasible solution of $\OPT( \bv+\mathcal{E} \textbf{1})$ and therefore
\begin{align*}
    \OPT( \bv+\mathcal{E} \textbf{1}) \leq \sum_{i=1}^k \iprod{\cv^{(i)}} {\xv^{(i)}_{*}} \leq \left(\sum_{i=1}^k \iprod{\cv^{(i)}} {\xv^{(i)}_{\textbf{L}}} \right) + \mathcal{E},
\end{align*}
where the last inequality follows from (\ref{eq:E}). Now, using (\ref{eq:wd_1}) with the above inequality, we obtain that
\begin{align*}
&\OPT(\bv+\mathcal{E} \textbf{1}) - \mathcal{E} \leq \DUAL(\bv) + \epsilon \\ 
\Rightarrow & \Delta \leq \mathcal{E} + \OPT(\bv) + \epsilon - \OPT(\bv+\mathcal{E} \textbf{1}),
\end{align*}
where the second inequality is obtained by subtracting $\OPT(\bv)$ on both sides.
Letting $\epsilon \to 0$ yields the desired statement.
\end{proof}
\end{theorem}

\section{Error bounds via Shapley-Folkman Lemma}\label{sec:newbnd}
In this section (and for the rest of the paper), we make an additional assumption:

\begin{asp}
    \label{asp:pointed}
    Let $A_i \subseteq \mathbb{R}^n$, $i \in [k]$ be the sets that are the summands in the Shapley--Folkman lemma. For all $i \in [k]$, we assume that $\conv(A_i)$ is closed and pointed.
\end{asp}
Assumption~\ref{asp:pointed} guarantees that each 
$\conv(A_i)$ has at least one extreme point. Moreover, since $\conv(A_i)$ is closed, any such extreme point must lie in the original set $A_i$. 
The pointedness assumption is subsequently 
used in Lemma~\ref{lem:face_decompose} in the next section to describe how the Minkowski sum of pointed closed convex 
sets decomposes in terms of their facial structure. This decomposition provides a geometric  proof of the Shapley--Folkman Lemma.


\subsection{Geometric proof of Shapley-Folkman Lemma}

In this section, we give a proof of Shapley-Folkman Lemma that utilizes the facial structure of the convex hull of each nonconvex set.

For completeness, we begin by presenting the standard definition of a face of a convex set and related notions.
\begin{definition}
Given a closed convex set $C$, we call a closed convex set $F \subseteq C$ a \textbf{face} of $C$ if for any line segment $[a,b] \subseteq C$ such that $(a,b) \cap F \neq \emptyset$, we have $[a,b] \in F$. 
The dimension of \blue{a} face $F$ denoted by $\dim(F)$ is the dimension of its affine hull.
\end{definition}
Note that a convex set $C$ is a face of itself. Moreover a face of a face of $C$ is also a face of $C$. Finally, faces of closed convex sets are closed convex sets \cite{hiriart2004fundamentals}.

We begin with some simple observations from  linear algebra and convex analysis. The first claim directly follows from a rank-nullity argument.
\begin{claim}
    \label{claim_rank}
    Let $L_i \in \R^n$, $i \in [k]$ be linear subspaces, such that $\sum_{i = 1}^k\dim(L_i) > n$. Then there exists $d_i \in L_i$ for $i \in [k]$, with not all \blue{$d_i$'s} zero, such that $\sum_{i = 1}^k d_i = 0.$
\end{claim}

\begin{claim}\label{claim:propface}
Let $C \subseteq \mathbb{R}^n$ be a closed convex set, and let $F$ be a nonempty face of $C$.
If $x \in \textup{relbd}(F)$, then $x$ is contained in a face of $C$ of
strictly smaller dimension than $F$.
\end{claim}
\begin{proof}
Since $x \in \textup{relbd}(F)$, there exists a hyperplane $H$ in the affine hull of $F$ such that (i) $F$ is contained in one of the half-spaces defined by $H$, (ii) $x \in H$, and (iii) $F \not \subseteq H$. Therefore, by definition of a face, $F\cap H$ is a face of $F$, and thus a face of $C$. Moreover, $(iii)$ implies that the dimension of $F$ is at least 1 more than the dimension of $F \cap H$.
\end{proof}

\begin{lemma}
    \label{lem:face_decompose}
    Let $A_1,\dots,A_{k}$ be nonempty sets of $\R^n$ with $k \geq n+1$ such that $\conv(A_i)$ is closed and pointed.
    Then for every $\yv \in \conv\left( \sum_{i=1}^k A_i\right)$, there exist faces $F_i$ of $\conv(A_i)$ such that
    \begin{align*}
        \yv \in \sum _{i=1}^kF_i \text{ and } \sum_{i=1}^k \dim(F_i) \leq n.
    \end{align*}
    \begin{proof}
    For any $\yv \in \conv\left(\sum_{i=1}^k A_i\right)$, consider $\yv = \sum_{i=1}^k \yv_i$ where each $\yv_i$ lies in the relative interior some face $F_i$ of $\conv(A_i)$ \blue{(if $\yv_i$ lies on the boundary of some face, we can choose a smaller face that includes $\yv_i$)}. 
    If $\sum_{i=1}^k \dim(F_i) > n$, we show how to find a different $\yv_i'$ from  $\conv(A_i)$ so that the total dimension $\sum_{i=1}^k \dim(F'_i)$, where  $F'_i$ is a face of $\conv(A_i)$ containing $\yv_i'$, 
    decreases while ensuring $\yv = \sum \yv_i'$. 
    Let $G_i$ be the linear subspace obtained by translating the affine hull of $F_i$ to the origin. If $\sum_{i=1}^k \dim(G_i) = \sum_{i=1}^k \dim(F_i) > n$, applying Claim \ref{claim_rank}, we can find $d_i \in G_i$, not all zero, such that  $ \dv_1 + \dv_2 + \cdots + \dv_k = \textbf{0}$.
Let $\yv'_i(\lambda) := \yv_i + \lambda d_i$ for $i \in [k]$. Since each $\conv(A_i)$ is pointed, there exists some choice of $\lambda$ such that for at least one $i* \in [k]$ we have that $\yv_{i*}'(\lambda)$ belongs to the relative boundary of $F_{i*}$ while keeping other $\yv_j'(\lambda)$ in $F_j$ for all $j \in [k]$. By Claim~\ref{claim:propface}, $\yv_{i*}'(\lambda) \in F'_{i*}$, a face of $\conv{(A_i)}$ such that $\textup{dim}(F'_{i*}) < \dim(F_{i*})$. Letting $F'_j  = F_j$ for all $j \not\eq i*$, the total dimension $\sum_{i=1}^k \dim(F'_i)$ decreases and we still ensure that $\yv = \sum_{i=1}^k \yv_i'$ for some $\yv'_i \in \conv(A_i)$.
\end{proof}

\end{lemma}

\begin{proof} \emph{of Shapley--Folkman Lemma using Lemma~\ref{lem:face_decompose}} For every $\yv \in \conv\left( \sum_{i=1}^k A_i\right)$, Lemma~\ref{lem:face_decompose} implies that there exist faces $F_i$ of $\conv(A_i)$ for $i \in [k]$, such that $ \yv \in \sum _{i=1}^kF_i \text{ and } \sum_{i=1}^k \dim(F_i) \leq n.$
Since $\dim(F_i) \in \{0,1,\dots,n\}$ and $\sum_{i=1}^k \dim(F_i) \leq n$, it implies that at least $k-n$ faces have dimension $0$. These dimension-$0$ faces are extreme points of corresponding $\conv(A_i)$s, therefore must lie in $A_i$. Let $\mathcal{I} = \{i \in[k]:\dim(F_{i}) > 0\}$ where from the above we have that $|\mathcal{I} | \leq n$. Moreover, we have,
\begin{align*}
    \yv \in \sum_{i \in \mathcal{I}} \conv(A_i) + \sum_{i \in [k] \setminus \mathcal{I} } A_i.
\end{align*}
\end{proof}

We note here that Lemma \ref{lem:face_decompose} can be used to strengthen the refined Shapley-Folkman Lemma recently proposed in \cite{bi2016refined}.


\subsection{Error bounds via randomized rounding}

Let $A_1, A_2, \dots, A_k$ be nonempty set of $\mathbb{R}^n$ such that $\conv(A_i)$ is closed and pointed. 
For any $\yv \in \conv\left(\sum_{i=1}^k A_i\right)$, Lemma~\ref{lem:face_decompose} implies the existence of faces $F_i$ of $\conv(A_i)$ such that $
    \yv = \sum_{i=1}^k \yv_i \text{ where } \yv_i \in F_i,~\forall i \in [k].$
By definition, $\Phi\left(\sum_{i=1}^k A_i\right)$ quantifies how close $\yv \in \sum_{i=1}^k \conv(A_i)$ is to some point $\yv' \in \sum_{i=1}^k A_i$ in the Euclidean norm, that is, $\|\yv - \yv'\|_2 \le \Phi(\sum_{i=1}^k A_i)$. 

The central idea of this section is that a suitable $\yv'$ can be constructed by randomly sampling points from each $F_i$. We first establish some preliminary results that facilitate this sampling approach.


\begin{claim}
    \label{claim_close_on_face}
    Let $C$ be a closed convex set and let $F$ be a face of $C$. \blue{For any $\xv \in F$, suppose $\xv$ is a convex combination of $\xv_i$s where $\xv_i \in C$, that is, $\xv = \sum_{i} \lambda_i \xv_i$,} with $\lambda_i \geq   0$ and $\sum_{i} \lambda_i = 1$. 
    \blue{Then $\xv_i \in F$ for every $i$ with $\lambda_i > 0$}
\begin{proof}
Directly follows from the definition of a face.
\end{proof}

\end{claim}



\begin{claim}
    \label{claim_variance_of_face}
    Let $A$ be a nonempty set of $\R^n$. If $\Xi(A) < \infty$, then
    for every face $F$ of $\conv(A)$, $\Xi(A \cap F) \leq \Xi(A)$.
\begin{proof}
\blue{For a nonempty $B \subseteq \R^n$ and $\xv \in \conv(B)$, write $\rho_B(\xv) := \inf\{\Rd(T) : T \subseteq B,\ \xv \in \conv(T)\}$, so that $\Xi(B) = \sup_{\xv \in \conv(B)} \rho_B(\xv).$ Fix $\xv \in \conv(A \cap F) \subseteq \conv(A) \cap F$, where the last containment follows from convexity of $F$. Let $\varepsilon > 0$. As $\rho_A(\xv) \leq \Xi(A)$, by definition of the infimum there is $T \subseteq A$ with $\xv \in \conv(T)$ and $\Rd(T) \leq \Xi(A) + \varepsilon$. Write $\xv = \sum_i \lambda_i \xv_i$ with $\xv_i \in T$, $\lambda_i \geq 0$ and $\sum_i \lambda_i = 1$, and set $T' := \{\xv_i : \lambda_i > 0\}$. Since $\xv \in F$ and $F$ is a face of $\conv(A)$, Claim~\ref{claim_close_on_face} gives $T' \subseteq F$, and hence $T' \subseteq A \cap F$. Moreover $\xv \in \conv(T')$ and, by monotonicity of $\Rd$, $\Rd(T') \leq \Rd(T) \leq \Xi(A) + \varepsilon$. In particular $T'$ is admissible in the infimum defining $\rho_{A \cap F}(\xv)$, so $  \rho_{A \cap F}(\xv) \leq \Rd(T') \leq \Xi(A) + \varepsilon.$ Since $\varepsilon > 0$ was arbitrary and $\Xi(A)$ does not depend on $\varepsilon$, we obtain $\rho_{A \cap F}(\xv) \leq \Xi(A)$. Taking the supremum over $\xv \in \conv(A \cap F)$ yields $\Xi(A \cap F) \leq \Xi(A)$.}
\end{proof}
\end{claim}

The above claim is interesting because, in general $\Xi(A\cap B)$ may not be smaller than $\Xi(A)$ or $\Xi(B)$. 

\begin{lemma}
\label{lem:sample_face} Let $A \subseteq \R^n$. Then for every $\xv \in \conv(A)$, there exists a finite distribution $\mathcal{D}$ such that $\mathcal{D}$ is supported on $A$, $\E_{Z \sim \mathcal{D}}(Z) = \xv$ and $\E(\norm{Z-\E(Z)}_2^2) \leq \Xi^2(A)$ .
\begin{proof}
    For any $\xv \in \conv(A)$, by the definition of $\Xi(A)$, there exist some $T \subseteq A$ such that $\xv \in \conv(T)$ and $\Rd(T) \leq \Xi(A)$. Since the statement of the Lemma is invariant under translation, we may translate $A$ so that $\norm{\yv}_2 \leq \Xi(A)$ for all $\yv \in T$.
    By Carathéodory's theorem, $\xv$ can be written as a convex combination of $\pv_i$ and each $\pv_i \in T$, that is $\xv = \sum_{i = 1}^{n+1}\lambda_i \pv_i \text{ where } \pv_i \in T$ and $\lambda$ belongs to the $n+1$ dimensional standard simplex.
    We define \blue{a} random variable $Z \sim \mathcal{D}$ that $Z = \pv_i$ with probability $\lambda_i$. It is clear that $\E_{Z \sim D}(Z) = \xv$ and 
    \begin{align*}
        \E(\norm{Z-\E(Z)}_2^2) = \E(\norm{Z}_2^2) - \norm{\E Z}_2^2 \leq  \E(\norm{Z}_2^2) \leq \Xi^2(A),
    \end{align*}
    where the last line uses that fact that  $Z$ only has support over $T$ and $\norm{\yv}_2 \leq \Xi(A)$ for all $\yv \in T$.
\end{proof}

\end{lemma}



\begin{lemma}
    \label{lem_SFL_quant_face}
    Let $A_1, A_2, \dots, A_k$ be nonempty subsets of $\mathbb{R}^n$. Let $F_i$ be a face of $\conv(A_i)$ for all $i \in [k]$ and let $\mathcal{I} := \{i \in [k] : \dim(F_i)>0\}$. Then it follows that
    \begin{align*}
        \Phi\left(\sum_{i=1}^k A_i \cap F_i  \right) \leq \sqrt{\sum_{i \in \mathcal{I}} \Xi^2(A_i \cap F_i)} \leq \sqrt{\sum_{i \in \mathcal{I}} \Xi^2(A_i)}.
    \end{align*}
    
\begin{proof}
    For every $\xv \in \conv\left( \sum_{i=1}^k A_i \cap F_i  \right)=\sum_{i=1}^k \conv(A_i \cap F_i)$ (the equality follows from Lemma~\ref{claim_convex_and_sum}), there \blue{exists} some $\xv_i$ such that $\xv = \sum_{i=1}^k \xv_i$ and $\xv_i \in \conv(A_i \cap F_i)$. We aim to construct $\yv_i$ such that $\yv_i \in A_i \cap F_i$ and $\norm{ \left( \sum_{i=1}^k \xv_i \right) - \left( \sum_{i=1}^k \yv_i \right) }_2 \leq \sqrt{\sum_{i \in \mathcal{I}} \Xi^2(A_i \cap F_i)}$. 
    
    For any $i \not\in \mathcal{I}$, $\dim(F_i) = 0$ and therefore $\xv_i$ is an extreme point of $\conv(A_i)$, so $\xv_i \in A_i \cap F_i$. 
    For $i \in \mathcal{I}$, by Lemma \ref{lem:sample_face}, there exists a distribution $Z_i \sim \mathcal{D}_i$ such that $\E_{Z_i\sim \mathcal{D}_i}(Z_i) = \xv_i$ and $\mathcal{D}_i$ is supported on $A_i \cap F_i$ and $\E(\norm{Z_i-\xv_i}_2^2) \leq  \Xi^2(A_i \cap F_i)$ 
    Then we (randomly) construct $\yv = \yv_1 + \yv_2+\dots \yv_k$ in the following way:
\begin{align*}
   \blue{\text{independently choose }} \yv_i = \begin{cases}
         \xv_i & \text{ if } i \in [k] \setminus \mathcal{I}, \\
         Z_i \sim \mathcal{D}_i & \text{ if } i \in \mathcal{I} 
    \end{cases}.
\end{align*}
Therefore, it follows that
\begin{align*}
    \E\left(\norm{\xv-\yv}_2^2\right) & = \E\left(\norm{\sum_{i \in \mathcal{I}}(  \xv_i-Z_i)}_2^2\right) \\
    & = \sum_{i \in \mathcal{I}} \E\left(\norm{(  \xv_i-Z_i)}_2^2\right) & \blue{\text{(because of independence)}} \\
    & \leq \sum_{i \in \mathcal{I}} \Xi^2(A_i \cap F_i).
\end{align*}
Therefore there exists a choice of $\yv$ such that 
$$\norm{ \left( \sum_{i=1}^k \xv_i \right) - \left( \sum_{i=1}^k \yv_i \right) }_2 \leq \sqrt{\sum_{i \in \mathcal{I}} \Xi^2(A_i \cap F_i)}.$$ Since $\xv$ is arbitrary and by Claim \ref{claim_variance_of_face}, this implies that $\Phi\left(\sum_{i=1}^k A_i \cap F_i  \right) \leq \sqrt{\sum_{i \in \mathcal{I}} \Xi^2(A_i \cap F_i)} \leq \sqrt{\sum_{i \in \mathcal{I}} \Xi^2(A_i)} .$
\end{proof}
\end{lemma}

Applying Lemma \ref{lem:face_decompose} and Lemma \ref{lem_SFL_quant_face} directly yields the following theorem.

\begin{theorem}\label{thm:generalrefinement}
    \label{thm_refine}
    Under conditions that each $\conv(A_i)$ is closed and pointed for all $i\in[k]$,
    we have that,
    $$\Phi\left(\sum_{i=1}^k A_i \right) \leq \sup\left\{\sqrt{\sum_{i =1}^k \Xi^2(A_i\cap F_i)}: F_i \text{ is a face of $A_i$ and $\sum_{i=1}^k \dim(F_i) \leq n$}\right\}.$$
    \blue{
    \begin{proof}
Fix any $
\yv \in \conv\left(\sum_{i=1}^k A_i\right).$
 By Lemma \ref{lem:face_decompose}, there exist faces $F_i$ of $\conv(A_i)$ such that
$
\yv \in \sum_{i=1}^k F_i\text{ and }
\sum_{i=1}^k \dim(F_i) \leq n.
$
Since $F_i$ is a face of $\conv(A_i)$, we have
$
F_i = \conv(A_i) \cap F_i = \conv(A_i \cap F_i),
$
and hence
$
\yv \in \sum_{i=1}^k \conv(A_i \cap F_i)
= \conv\left(\sum_{i=1}^k (A_i \cap F_i)\right),
$
where the last equality follows from Lemma \ref{claim_convex_and_sum}.
Applying Lemma~\ref{lem_SFL_quant_face} to the sets $A_i \cap F_i$, there exists some point $\zv$ in $\sum_{i=1}^k (A_i \cap F_i) \subseteq \sum_{i=1}^k A_i$ such that $\norm{\zv-\yv}_2 \leq \sqrt{\sum_{i \in \mathcal{I}} \Xi^2(A_i \cap F_i)}$ where $\mathcal{I} := \{i \in [k] : \dim(F_i)>0\}$.
Taking the supremum over all
$\yv \in \conv(\sum_{i=1}^k A_i)$ and over all possible faces $F_i$ satisfying
$\sum_{i=1}^k \dim(F_i)\leq n$ gives
\[
\Phi\left(\sum_{i=1}^k A_i\right)
\leq
\sup\left\{
\sqrt{\sum_{i=1}^k \Xi^2(A_i \cap F_i)}
:
F_i \text{ is a face of } \conv(A_i),
\ \sum_{i=1}^k \dim(F_i)\leq n
\right\}.
\]
\end{proof}
}
\end{theorem}

Theorem \ref{thm_refine} is a refinement on the classical error bounds,
due to the fact that $\Xi(A_i \cap F_i)$ can be significantly smaller than $\Xi(A_i)$ unless $F_i$ is a high-dimensional face. However, due to the global dimension-sum constraints, the number of such high-dimensional faces is limited.

A more concrete illustration is as follows. Consider any full-dimensional convex body $C \in \mathbb{R}^n$ with $\Rd(C)=R$, and \blue{define for all $i \in [k]$ the set $
A_i = C \setminus \textup{int}(C)$.}
In this case, every face $F_i$ of $\conv(A_i)$ that is not $\conv(A_i)$ itself, satisfies $
\Xi(A_i \cap F_i) = 0$ because each $F_i \cap A_i$ is a convex set. 
\blue{Therefore, $\Xi(A_i \cap F_i) \neq 0$, only when $F_i = \conv(A_i)$ {and in this case $\dim(F_i) = n$}.}
Now applying Theorem \ref{thm_refine} implies that
$\Phi\!\left(\sum_{i=1}^k A_i \right) \leq \Xi(A_1) \leq R,$
which is significantly smaller than the $\sqrt{n}R$ bound predicted by the classical error estimates derived from the Shapley--Folkman lemma.



\subsection{Improved error bounds via local geometry}

In this subsection, we seek sharper bounds on 
$\Phi\!\left(\sum_{i=1}^k A_i\right)$ by leveraging the geometric 
structure of the summands $A_i$. Our goal is to identify sufficient 
conditions under which
\[
\Phi\!\left(\sum_{i=1}^k A_i\right) \;\longrightarrow\; 0 
\quad \text{as } k \to \infty,
\]
that is, the Minkowski sum becomes asymptotically convex. 

For each $\yv = \sum_{i=1}^k \yv_i \in \conv\!\left(\sum_{i=1}^k A_i\right)$, Lemma~\ref{lem:face_decompose} implies that each $\yv_i$ lies in the Minkowski sum of certain faces $F_i$ of $\conv(A_i)$. 
The error bound presented in Lemma~\ref{lem_SFL_quant_face} can be interpreted as follows: we attempt to round $\yv$ as a whole to a point in $\sum_{i=1}^k (A_i \cap F_i)$. 
This perspective naturally suggests a potential further refinement of Lemma~\ref{lem_SFL_quant_face}, wherein we aim to round $\yv$ directly to a point in $\sum_{i=1}^k A_i$. 
However, this 
appealing idea encounters a significant challenge, as it requires full characterization of the nonconvex sets $A_i$, which is typically unavailable in practice. 

The key insight of our main result is the following: We do not need to understand the global structure of the nonconvex sets $A_i$. It is sufficient to instead focus on exploiting only the local geometric properties of each $A_i$.   


We next introduce a definition that plays a key role in identifying the required local geometry.

\begin{definition}
Let $A$ be a nonempty set. For \blue{a subset $Q \subseteq A$}, we call $H \subseteq A$ is \textit{hidden convex component} associated with $Q$ if
    \begin{enumerate}
        \item $H \subseteq A$ is a closed convex set;
        \item $Q \subseteq H$;
    \end{enumerate}
\end{definition}

\begin{figure}[thbp]
    \centering
    \includegraphics[width=0.5\linewidth]{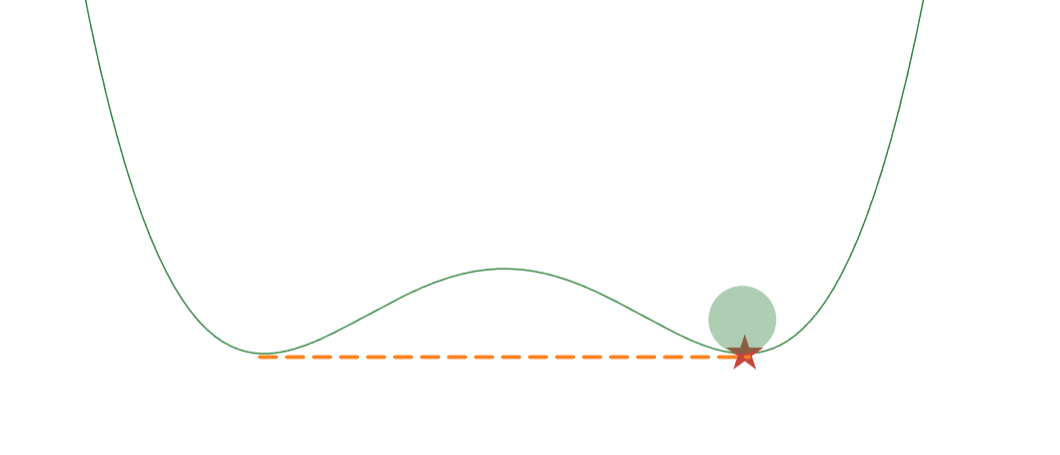}
    \caption{$A$ is the epigraph of the function $y=x^4-x^2$. The star point is the set $Q$. The $l_2$ ball contained in the epigraph is a hidden convex component.}
    \label{fig:hidden_convex}
\end{figure}
The definition above is illustrated via an example in Figure \ref{fig:hidden_convex}.

Such hidden convex component may not 
exist for general $Q$. On the other hand,
for each extreme point $\mathbf{v}_i$ of $\conv(A_i)$, there always exists a trivial hidden convex 
component---namely the singleton set $H_i = \{\mathbf{v}_i\}$. It is immediate that the existence 
of such trivial hidden convex components does not drive 
$\Phi\!\left(\sum_{i=1}^k A_i\right)$ to zero. For instance, take 
$A_i = \{0,1\}^n$ for all $i \in [k]$. Their Minkowski sum 
$\sum_{i=1}^k A_i$ is contained in the integer lattice $\{0,1,\dots,k\}^n$, and clearly 
$\Phi\!\left(\sum_{i=1}^k A_i\right)$ does not vanish.
A natural next attempt is to require each hidden convex component $H_i$ to have nonempty interior. 
However, as  discussed in Remark~\ref{rem:p_ball} the mere existence of hidden convex components with non-vanishing interior still does not guarantee any improvement in the convexification of the Minkowski sum. 
In fact, diminishing the nonconvexity measure 
$\Phi\!\left(\sum_{i=1}^k A_i\right)$ requires an appropriate smoothness condition of the hidden convex component $H_i$s.

\begin{definition}\cite{ben2001lectures}
A differentiable function $f(\cdot):\R^n \to \R$ is called $L$-smooth if
$f(\yv) \leq f(\xv) + \iprod{\yv-\xv}{\nabla f(\xv)} + \dfrac{L}{2} \norm{\yv-\xv}_2^2,\forall \yv \in \R^n$, where $\nabla f(\cdot)$ is the gradient.
\end{definition}

\begin{definition}\cite{hiriart2004fundamentals}
Let $C \subseteq \mathbb{R}^n$ be a nonempty closed convex set with $\text{0}$ in its interior.
The \emph{gauge function} associated with $C$ is
\[
    \norm{\xv}_{C}
    := 
    \inf\{\lambda > 0 : x \in \lambda C\},
    \qquad x \in \mathbb{R}^n.
\]
\end{definition}

For a hidden convex component $H$ with $\textbf{0}$ in its interior, we quantify its smoothness by the smoothness of the square gauge function $\norm{\cdot}_{H}^2$, which is related to the concept of $(2,D)$-smooth in functional analysis \cite{molinaro2023strong}. 


To utilize the local geometry, we first overload the definition of $\Phi^{K}(\cdot)$: given a closed convex set $K$ with $\textbf{0}$ in its interior, we let
    \begin{align*}
        \Phi^{K}(A,D) := \sup_{\xv \in A} \inf_{\yv \in D} \norm{\xv-\yv}_{K}.
    \end{align*}
It is therefore straightforward to see that $\Phi^{K}(A) = \Phi^{K}(\conv(A),A)$. Moreover, we note that $\Phi^{K}(A,D)$ is invariant under translation of $A$ and $D$ by the same vector.

\begin{lemma}
    \label{lem:local_geom}
    Let $A_1, A_2, \dots, A_k$ be nonempty subsets of $\mathbb{R}^n$ and $\Xi(A_i) \leq \beta,\forall i \in [k]$. Let $F_i$ be a face of $\conv(A_i)$, for all $i \in [k]$ and let $\mathcal{I} := \{i \in [k] : \dim(F_i)>0\}$. Let $H_i$ be a hidden convex component associated with $F_i$ for all $i \in [k] \setminus \mathcal{I}$ and $H = \sum\limits_{i \in [k] \setminus \mathcal{I}} H_i$. Suppose $H$ has nonempty interior, let $\mathcal{H}$ be a translation of $H$ that contains $\textbf{0}$ in its interior. If we have $\norm{\cdot}_{\mathcal{H}}^2$ is $L$-smooth, then 
 \begin{align*}
    \Phi^{\mathcal{H}}\left(\sum_{i=1}^k F_i,\sum_{i=1}^k A_i\right) \leq \sqrt{1 +  \dfrac{L}{2} \sum_{i \in \mathcal{I}} \Xi^2(A_i \cap F_i)}-1. 
    \end{align*}
    \begin{proof}
$\Phi^{\mathcal{H}}\left(\sum_{i=1}^k F_i,\sum_{i=1}^k A_i\right)$ is equivalent to the maximum possible distance (under $\norm{\cdot}_{\mathcal{H}}$) from points in $\sum_{i=1}^k F_i$ to $\sum_{i=1}^k A_i$. If $\sum_{i=1}^k F_i \subseteq \mathcal{F}$ for some $\mathcal{F} \subseteq \R^n$, it is clear that $\Phi^{\mathcal{H}}\left(\sum_{i=1}^k F_i,\sum_{i=1}^k A_i\right) \leq \Phi^{\mathcal{H}}\left(\mathcal{F},\sum_{i=1}^k A_i\right)$.
We derive an upper bound of this quantity by choosing $\mathcal{F} := \left(\sum\limits_{i \in \mathcal{I}} F_i \right) + \left( \sum\limits_{i \in [k] \setminus \mathcal{I}} H_i \right)$. 
We will show that for any point in $\left(\sum\limits_{i \in \mathcal{I}} F_i \right) + \left( \sum\limits_{i \in [k] \setminus \mathcal{I}} H_i \right)$ there exists a point in $\left(\sum\limits_{i \in \mathcal{I}} A_i \cap F_i \right) + \left( \sum\limits_{i \in [k] \setminus \mathcal{I}} H_i \right) \subseteq \sum_{i=1}^k A_i$ within the desired distance.  

Since the target statement is invariant under equal translation of $\sum_{i=1}^nF_i$ and $\sum_{i = 1}^nA_i$, we first translate $A_i$ so that $\sum\limits_{i \in [k] \setminus \mathcal{I}} H_i = \mathcal{H}$. Therefore, any point $\pv$ in $\left(\sum\limits_{i \in \mathcal{I}} F_i \right) + \left( \sum\limits_{i \in [k] \setminus \mathcal{I}} H_i \right)$, can be written as
$\pv = \left(\sum\limits_{i \in \mathcal{I}} \pv_i\right) + \vv$
where $\pv_i \in F_i$ and $\vv \in \mathcal{H}$. By Lemma \ref{lem:sample_face}, we can equip each $F_i$ with a distribution $Z_i \sim D_i$ with support on $A_i \cap F_i$ such that $\E_{Z_i \sim \mathcal{D}_i}(Z_i) = \pv_i$ and $\E(\norm{Z_i-\E(Z_i)}_2^2) \leq \Xi^2(A_i \cap F_i)$.
\blue{Fix any point $\pv$, we complete the proof in two steps:
\begin{enumerate}
    \item We first find some point $\sum_{i \in \mathcal{I}} Z_i$ via randomized rounding that is close to $\pv$.
    \item We then further move $\sum_{i \in \mathcal{I}} Z_i$ to a point in
    $\left(\sum_{i \in \mathcal{I}} (A_i \cap F_i)\right)
    + \left(\sum_{i \in [k] \setminus \mathcal{I}} H_i\right)$
    that is closer to $\pv$.
\end{enumerate}
}
\blue{For the first step, we first try to bound the expected distance:}
\begin{align*}
 \E\left( \norm{\pv-\left(\sum_{i \in \mathcal{I}} Z_i\right)}_{ \mathcal{H}}^2 \right) & = \E\left( \norm{\vv-\left(\sum_{i \in \mathcal{I}}Z_i-\pv_i\right)}_{ \mathcal{H}}^2 \right).
\end{align*}
Let $\gv:=\nabla \norm{\vv}^2_{\mathcal{H}}$ and applying the definition of $L$-smoothness, this further yields that
\begin{align*}
    \E\left( \norm{\vv-\left(\sum_{i \in \mathcal{I}} Z_i-\pv_i\right)}_{\mathcal{H}}^2 \right) & \leq \E\left( \norm{\vv}_{\mathcal{H}}^2 \right) + \E\left(\iprod{-\sum_{i \in \mathcal{I}} \left(Z_i - \pv_i\right)}{\gv}\right)+\dfrac{L}{2} \E\left( \norm{\left(\sum_{i \in \mathcal{I}}Z_i-\pv_i\right)}^2_{2} \right)  \\
    & = \E\left( \norm{\vv}_{\mathcal{H}}^2 \right) +  \frac{L}{2}\E\left( \norm{\left(\sum_{i \in \mathcal{I}}Z_i-\pv_i\right)}_2^2 \right)  \\
    & = \E\left( \norm{\vv}_{\mathcal{H}}^2 \right) + \frac{L}{2} \sum_{i \in \mathcal{I}} \E(\norm{Z_i-\pv_i}_2^2) \ \ \blue{\text{(because of independence)}} \\
    & \leq 1 +  \frac{L}{2} \sum_{i \in \mathcal{I}} \Xi^2(A_i \cap F_i).
\end{align*}
\blue{
This implies there exists some point $\zv$ in $\sum_{i \in \mathcal{I}} A_i \cap F_i$ whose distance defined by $\norm{\cdot}_{\mathcal{H}}$ to $\pv$ is at most $\sqrt{1 +  \frac{L}{2} \sum\limits_{i \in \mathcal{I}} \Xi^2(A_i \cap F_i)}$. That is, $\pv-\zv= \norm{\pv-\zv}_{\mathcal{H}} \rv$ for some $\rv \in  \mathcal{H}$ where $\norm{\pv-\zv}_{\mathcal{H}}\leq\sqrt{1 +  \frac{L}{2} \sum\limits_{i \in \mathcal{I}} \Xi^2(A_i \cap F_i)}$. Finally, we can choose $\zv+\rv \in \left(\sum_{i \in \mathcal{I}} A_i \cap F_i\right)+\mathcal{H}$ so that $\pv-(\zv+\rv )=(\norm{\pv-\zv}_{\mathcal{H}} -1) \rv$. This implies that
the $\norm{\cdot}_{\mathcal{H}}$ distance of $\pv$ to $\left(\sum_{i \in \mathcal{I}} A_i \cap F_i\right)+\mathcal{H}$ is at most $\sqrt{1 +  \frac{L}{2} \sum\limits_{i \in \mathcal{I}} \Xi^2(A_i \cap F_i)}-1$.}
\end{proof}
\end{lemma}

\begin{theorem}
    \label{thm_rem_convex_ball}
    Let $A_1, A_2, \dots, A_k$ be nonempty subsets of $\mathbb{R}^n$ and $\Xi(A_i) \leq \beta$, for all $i \in [k]$. \blue{Further assume that that $\conv(A_i)$ is closed and pointed}. Suppose every extreme point of $A_i$ has an associated $l_2$--ball with radius $r_i$ ($r_i$ can be zero) as a hidden convex component. Let $r_* = \min_{\mathcal{I}\subseteq [k]:|\mathcal{I}|=n} \sum_{i \in [k] \setminus \mathcal{I}} r_i$. Then it follows that
    \begin{align*}
        \Phi\left(\sum_{i=1}^k A_i\right) \leq  \sqrt{r_*^2 +  n \beta^2 }- r_*.
    \end{align*}
If there exists some $r > 0$ such that $r_i \geq r$ for all $i \in [k]$, then 
\begin{align*}
    \Phi\left(\sum_{i=1}^k A_i\right) \leq  \sqrt{ (k-n)^2r^2 + n \beta^2 }- (k-n)r.
\end{align*}
As $k\to\infty$, it is clear the above $\Phi\left(\sum_{i=1}^k A_i\right) \to 0$.
\begin{proof}
To bound $\Phi\left(\sum_{i=1}^k A_i\right)$, it suffices to show that
    for any $\pv \in \sum_{i=1}^k \conv(A_i)$, we have $\Phi\left(\pv,\sum_{i=1}^k A_i\right) \leq \sqrt{r_*^2 +  \dfrac{1}{2} n \beta^2 }- r_*$. By Lemma \ref{lem:face_decompose}, there exists face $F_i$ of $\conv(A_{i})$ such that $\pv \in \sum_{i=1}^k F_i$ and $\sum_{i=1}^k \dim(F_i) \leq n$. 
 
 By the hypothesis, after suitable translation of all the sets, we can find $H_i = \mathcal{B}(\textbf{0},r_i)$ as hidden convex component associated with $F_i$ for all $i \in [k] \setminus \mathcal{I}$ such that $\mathcal{H} = \sum\limits_{i \in [k] \setminus \mathcal{I}} H_i = \mathcal{B}(\textbf{0},r_*)$. 
    Then $\norm{\xv}_{\mathcal{H}} = \dfrac{1}{r_*} \norm{\xv}$ and $\norm{\cdot}_\mathcal{H}^2$ is $\dfrac{\blue{2}}{r_*^2}$-smooth. Lemma \ref{lem:local_geom} implies that
     \begin{align*}
         \Phi^{\mathcal{H}}\left(\sum_{i=1}^k F_i,\sum_{i=1}^k A_i\right) \leq \sqrt{1 +  \dfrac{1}{r_*^2} \sum_{i \in \mathcal{I}} \Xi^2(A_i \cap F_i)}-1.
     \end{align*}
This implies that
\begin{align*}
    \Phi\left(\sum_{i=1}^k F_i,\sum_{i=1}^k A_i\right) & \leq \sqrt{r_*^2 + \sum_{i \in \mathcal{I}} \Xi^2(A_i \cap F_i)}- r_* \\
    & \leq \sqrt{r_*^2 + \sum_{i \in \mathcal{I}} \beta^2}- r_*,
\end{align*}
where the last inequality uses Claim \ref{claim_variance_of_face}. Since $\pv$ is arbitrary, this implies that $ \Phi\left(\sum_{i=1}^k A_i\right) \leq  \sqrt{r_*^2 + n \beta^2 }- r_*$. The rest of statement directly follows that $r_* \geq (k-n)r$.
\end{proof}
\end{theorem}

\blue{Note that in the above result, ideally the set $\mathcal{I}$ used in $r_*$ would be the index set of components whose associated face has positive dimension, that is, $\mathcal{I} := \{i \in [k] : \dim(F_i) > 0\}$. This set is not known a priori and depends on the point under consideration, so in the statement of the theorem we choose $\mathcal{I}$ according to the worst case.}



We end the section with several examples and remarks.

\begin{remark}
Lemma~\ref{lem:local_geom} can be generalized by quantifying the smoothness of the hidden convex components through the smoothness (possibly with respect to a different norm) of the $p$-th power of a gauge function $\|\cdot\|_{H}^{p}$ for some $p>1$, following the same probabilistic argument. 
We do not pursue this extension here, as it naturally requires alternative notions of nonconvexity beyond those considered in this work. 
For simplicity, we therefore focus on the classical nonconvexity measures $\Phi(\cdot)$ and $\Xi(\cdot)$ that appear in the literature on the Shapley--Folkman lemma, due to their direct relevance in estimating the duality gap~$\Delta$.

\end{remark}

\begin{remark}
\label{rem:p_ball}
In general, the smoothness of \blue{the} hidden convex component is required for establishing the vanishing of nonconvexity. To illustrate this phenomenon, consider the simple example $
A_0 + \sum_{i=1}^k \mathcal{B}_p(\mathbf{0},1)=A_0+\mathcal{B}_p(\mathbf{0},k)$
where $A_0 := \{(0,0),(0,1)\} \subseteq \mathbb{R}^2$ is a nonconvex set and $\mathcal{B}_p(\mathbf{0},1)$ denotes the unit $p$-norm ball. Figures~\ref{fig:ball_1}, \ref{fig:ball_1.5} and \ref{fig:ball_2} depict the resulting Minkowski sums for $k=1,3,$ and $6$, respectively.
We observe that the ``nonconvex hollow'' \blue{(that is, $\conv\left(A_0+\mathcal{B}_p(\mathbf{0},k)\right)\setminus A_0+\mathcal{B}_p(\mathbf{0},k)$)} in the set gradually vanishes as the number of summands increases when $p=1.5,2$. In contrast, this hollow persists for all $k$ in the case $p=1$. The underlying reason is that $\norm{\cdot}_{\mathcal{B}_{1}(\textbf{0},1)}^f$ is not smooth for any choice of $f > 1$.
\begin{figure}[!htb]
\minipage{0.32\textwidth}
  \includegraphics[width=\linewidth]{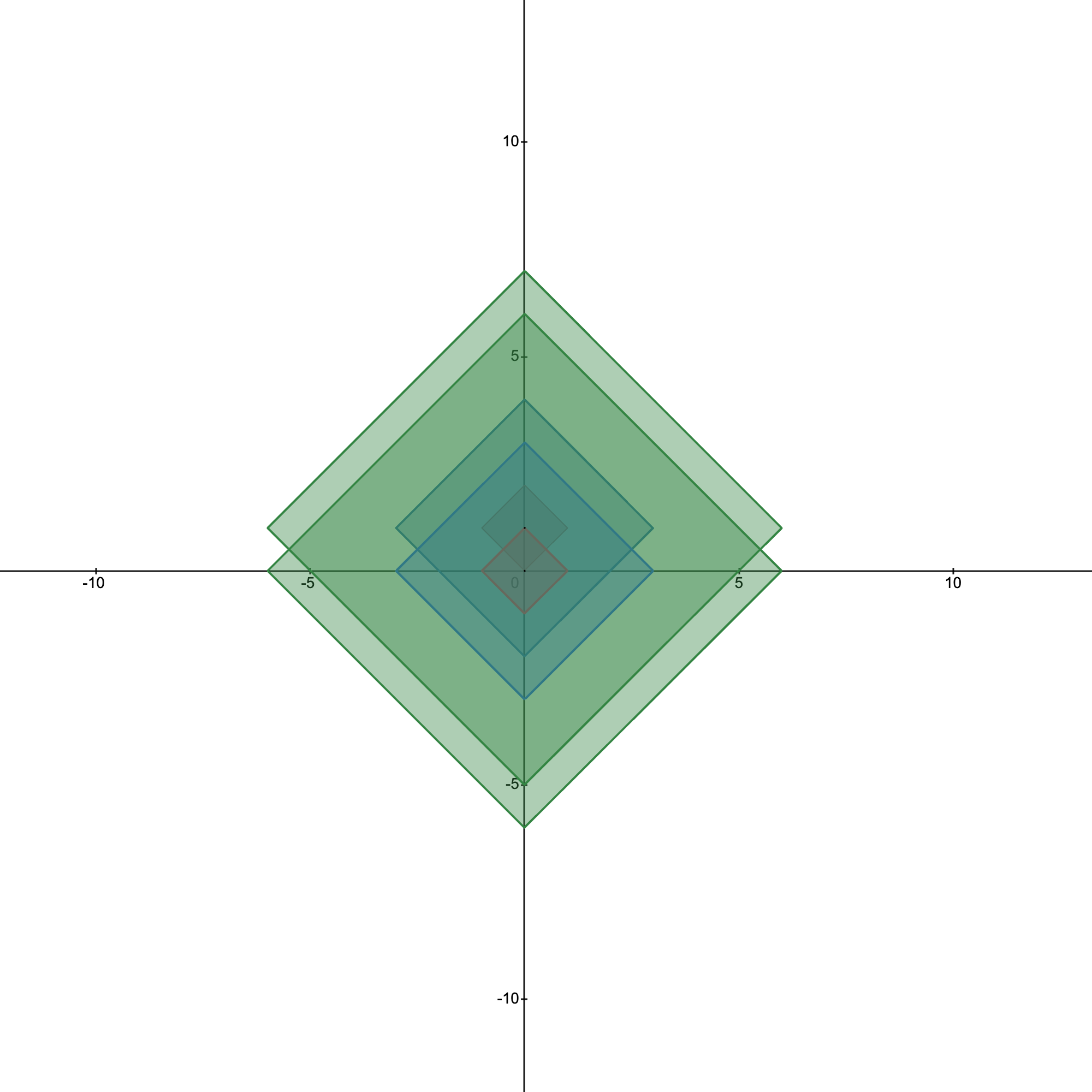}
  \caption{$p=1$}\label{fig:ball_1}
\endminipage\hfill
\minipage{0.32\textwidth}
  \includegraphics[width=\linewidth]{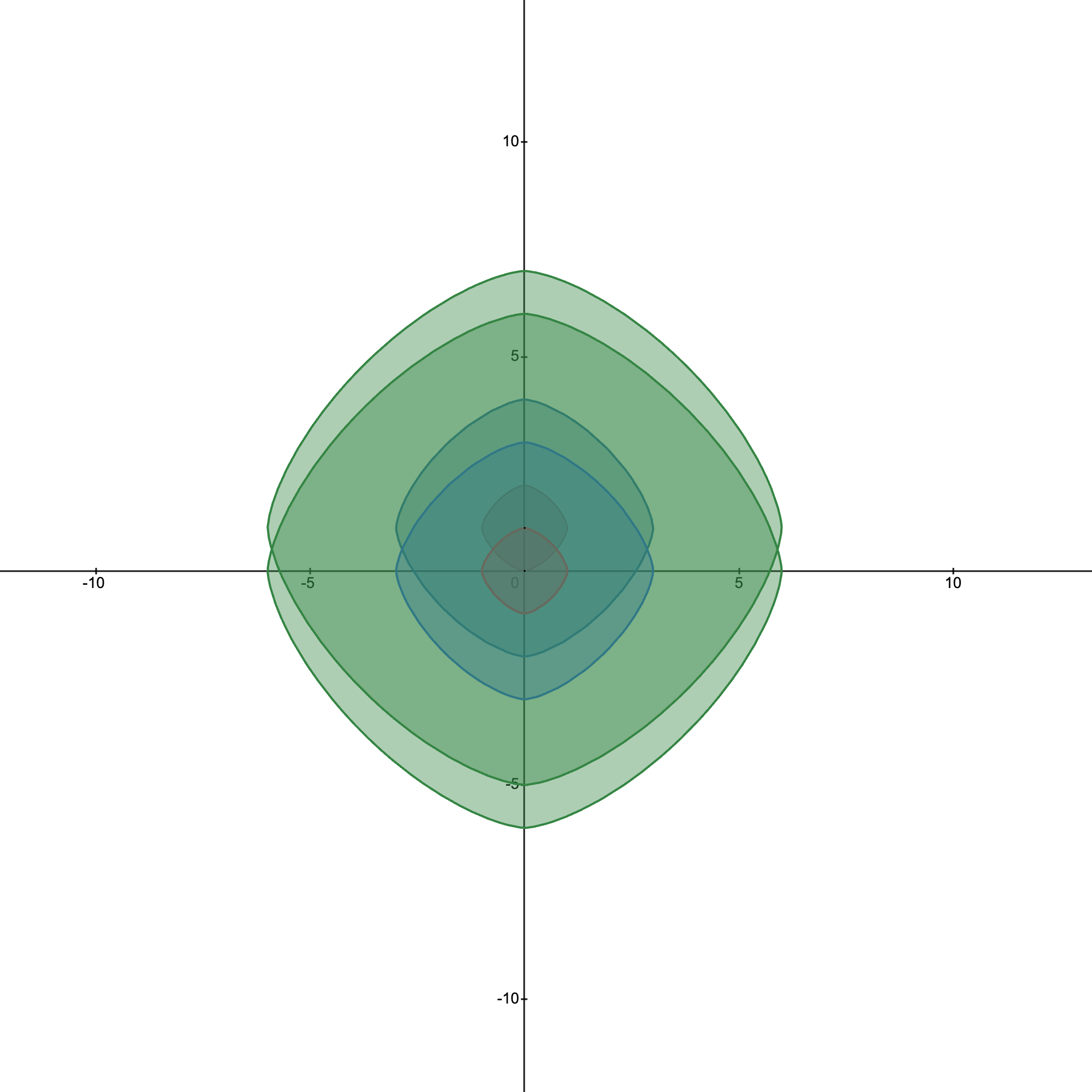}
  \caption{$p=1.5$}\label{fig:ball_1.5}
\endminipage\hfill
\minipage{0.32\textwidth}%
  \includegraphics[width=\linewidth]{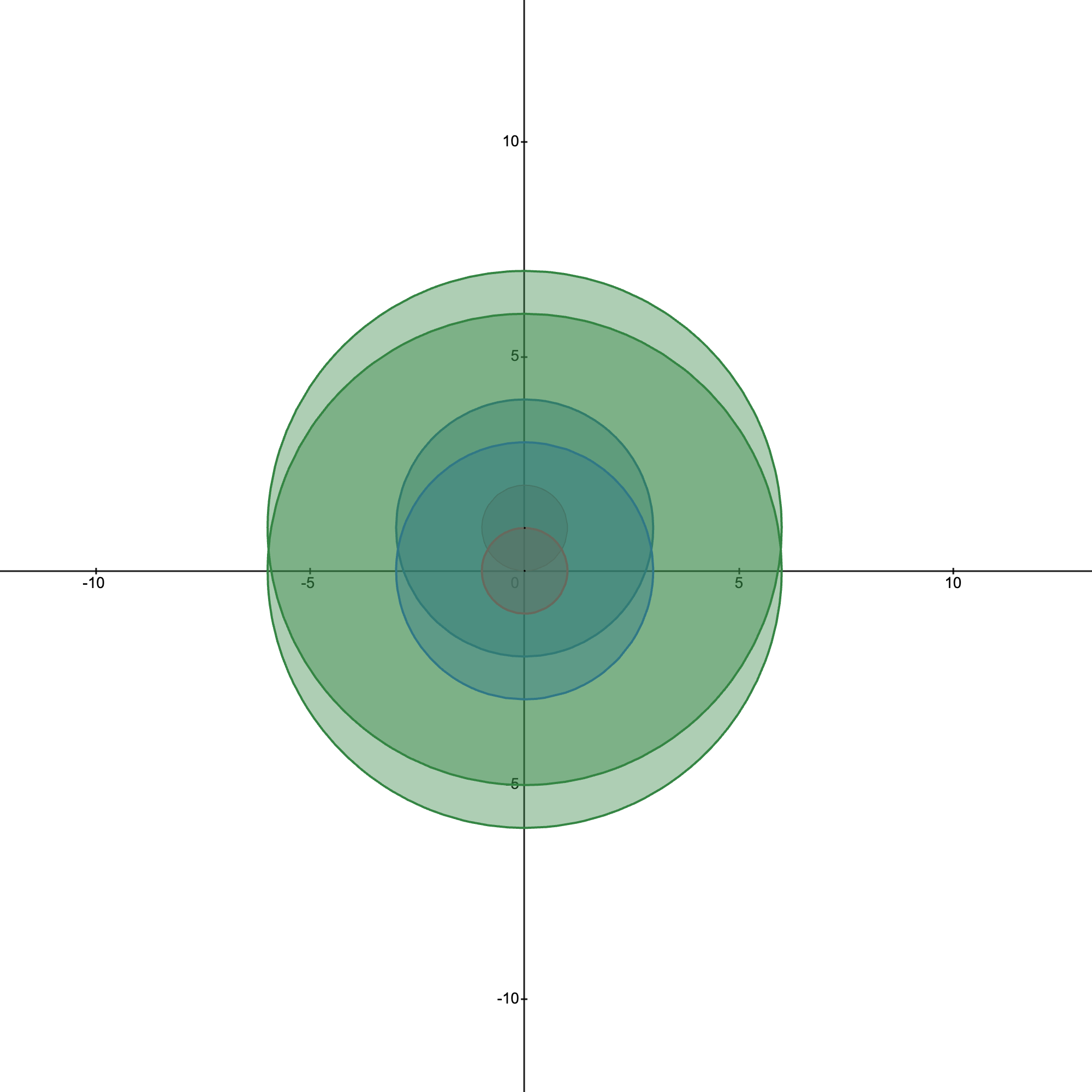}
  \caption{$p=2$}\label{fig:ball_2}
\endminipage
\end{figure}
\end{remark}

\begin{remark}
To ensure that the Minkowski sums converge to a convex set, it is sufficient that smooth, hidden convex components exist only in neighborhoods of the extreme points of \blue{the} convex hull of each nonconvex set. Such components need not exist around non-extreme points of the convex hulls.
To illustrate this, consider the $\mathbb{R}^2$ examples $T_1$ and $T_2$ shown in Figures~\ref{fig:desmos_1}, ~\ref{fig:desmos_2}, ~\ref{fig:desmos_3} and~\ref{fig:desmos_4}. Both sets contain some points on which no smooth hidden convex components are associated with. Nevertheless, one can verify that $\sum_{i=1}^k T_1$ converges to a convex set as $k \to \infty$, whereas $\sum_{i=1}^k T_2$ does not.
The difference arises from the behavior at extreme points: every extreme point of $\conv(T_1)$ admits a smooth hidden convex component (e.g., a ball contained in $T_1$ in a neighborhood of that point), while $\conv(T_2)$ possesses five extreme points for which no such smooth hidden convex component exists.
\begin{figure}[!htb]
\minipage{0.25\textwidth}
  \includegraphics[width=\linewidth]{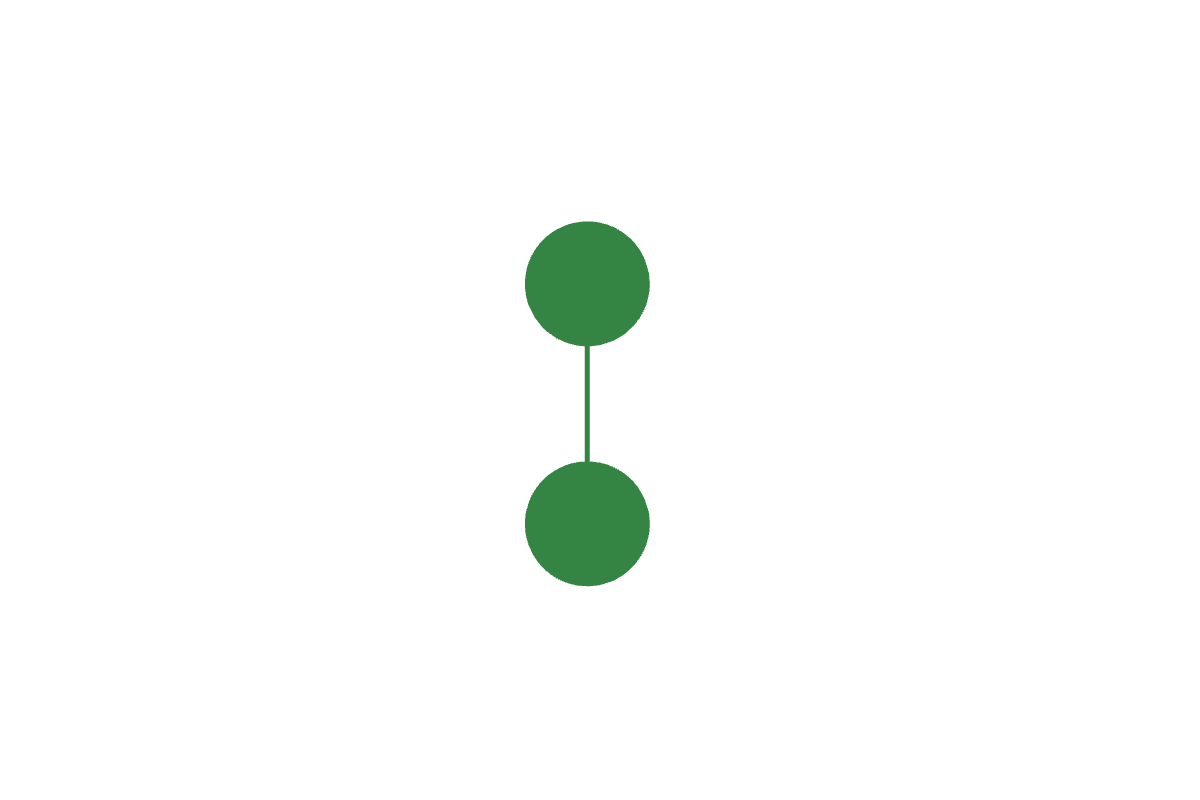}
  \caption{$T_1$}\label{fig:desmos_1}
\endminipage\hfill
\minipage{0.25\textwidth}
  \includegraphics[width=\linewidth]{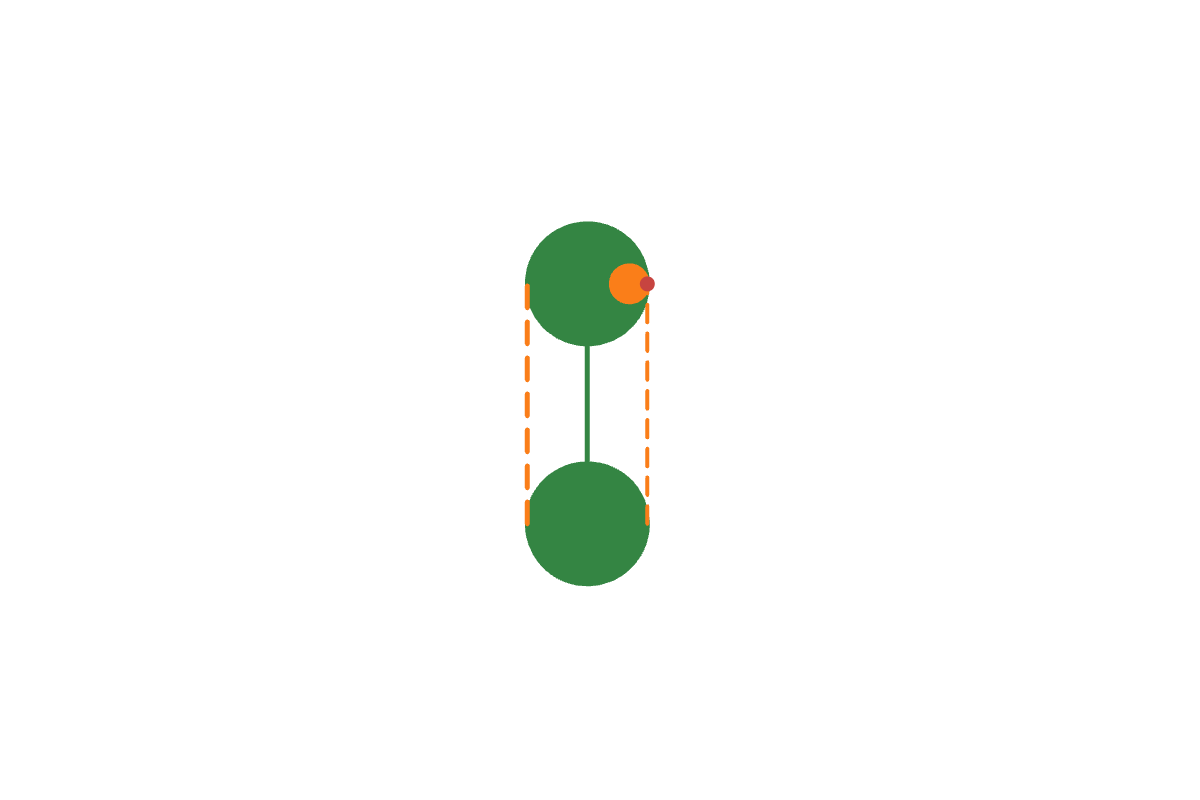}
  \caption{$\conv(T_1)$}\label{fig:desmos_2}
\endminipage\hfill
\minipage{0.25\textwidth}%
  \includegraphics[width=\linewidth]{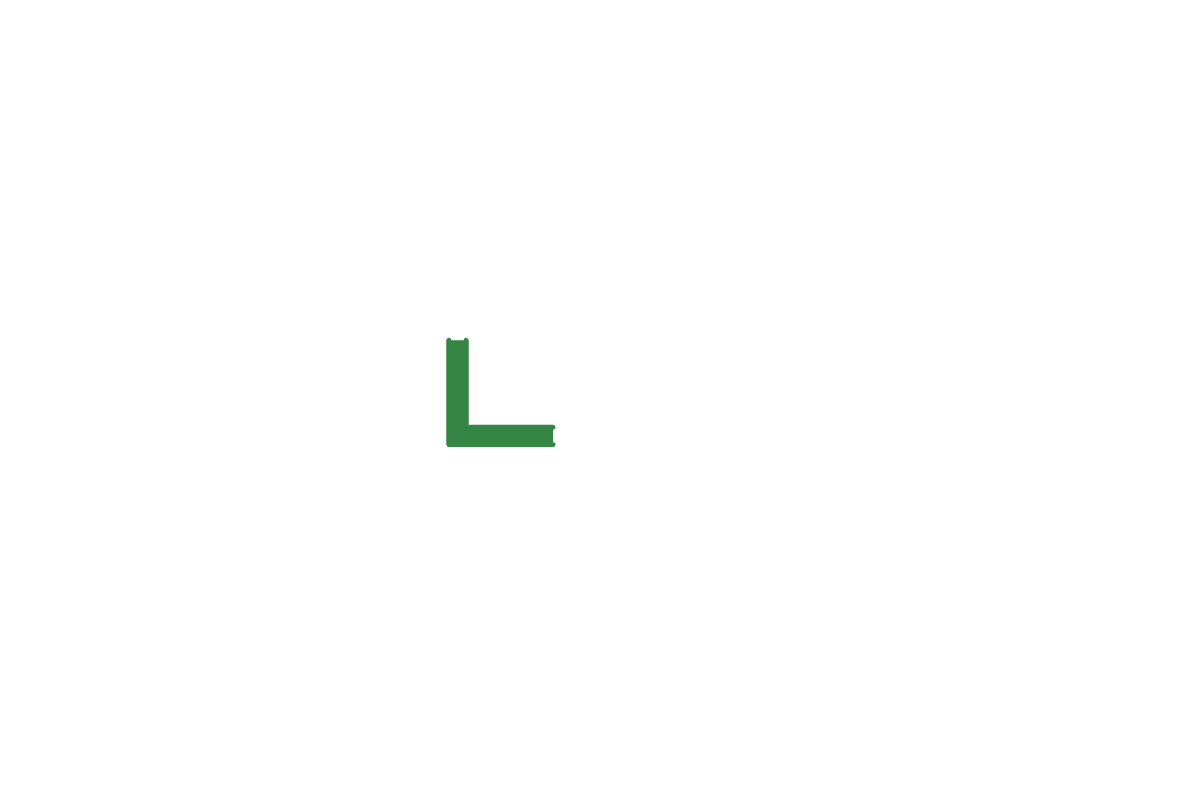}
  \caption{$T_2$}\label{fig:desmos_3}
\endminipage
\minipage{0.25\textwidth}%
  \includegraphics[width=\linewidth]{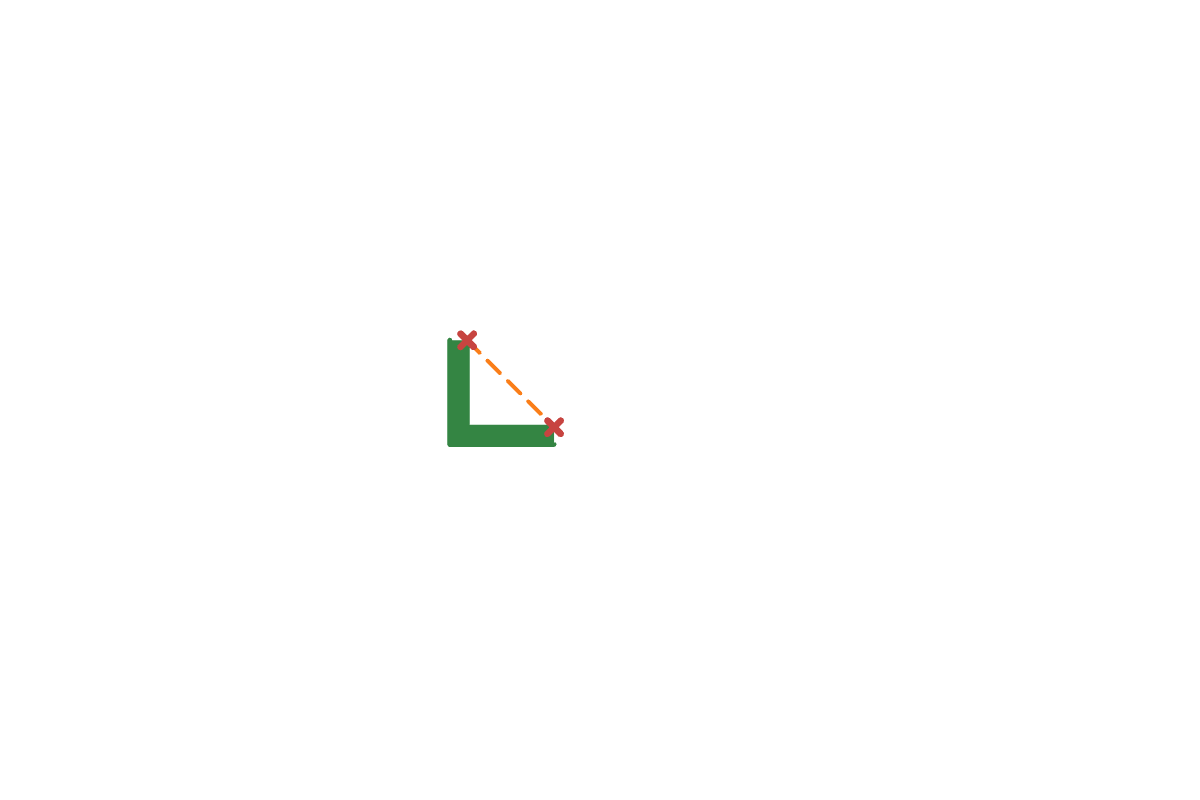}
  \caption{$\conv(T_2)$}\label{fig:desmos_4}
\endminipage
\end{figure}
\end{remark}

\section{Asymptotically tight Lagrangian dual of smooth nonconvex problems}\label{sec:app}

In this section, we aim to leverage Theorem~\ref{prob_smooth_sparsity} and Lemma~\ref{lem:local_geom} to establish the {asymptotic tightness} of problem~(\ref{prob_smooth_sparsity}). The key insight is that the smooth hidden convex components naturally emerge at the extreme points of the convex hull of the epigraphs of smooth functions. Under mild assumptions on the constraint matrices $B^{(i)}$, the smoothness of these hidden convex components is approximately preserved under projection, thereby implying that $\Phi(\mathcal{P})$ vanishes as $k$ increases. We rewrite (\ref{prob_smooth_sparsity}) as:
\begin{equation}
    \label{eq:sqarsity_constrainted}
    \begin{aligned}
         \OPT_{\textbf{s}}(\bv) := \inf & \sum_{i=1}^{k} t^{(i)} \\
    \text{s.t. } & \sum_{i=1}^k B^{(i)} \xv^{(i)} \leq \bv, \\
    & \norm{\xv^{(i)}}_0 \leq s_i,t_i \geq f^{(i)}(\xv^{(i)}),\forall i \in[k]. \\
    \end{aligned}
\end{equation}
In this case, we have 
\begin{align*}
   \mathcal{X}^{(i)}_{\textbf{s}} & := \left\{ (t,\xv) \in \R^{m+1} \;\middle\vert\;
   \begin{array}{@{}l@{}} \norm{\xv^{(i)}}_0 \leq s_i, t^{(i)} \geq f^{(i)}(\xv^{(i)}) \\
   \end{array}
   \right\}
\end{align*}
and the projected set $\mathcal{P}^{(i)}$ for (\ref{eq:sqarsity_constrainted}) is
\begin{align*}
    \mathcal{P}^{(i)}_{\textbf{s}}& := \left\{ (t,\dv) \in \R^{m+1} \;\middle\vert\;
   \begin{array}{@{}l@{}} \exists (t^{(i)},\xv^{(i)}) \text{ that } \norm{\xv^{(i)}}_0 \leq s_i, t^{(i)} \geq f^{(i)}(\xv^{(i)}) \\
   t = t^{(i)}, \dv = B^{(i)} \xv^{(i)}
   \end{array}
   \right\} \\
   & = \left\{ (t,\dv) \in  \R^{m+1} \;\middle\vert\;
   \begin{array}{@{}l@{}} \exists \xv^{(i)} \text{ that } \norm{\xv^{(i)}}_0 \leq s_i, \\
   t \geq f^{(i)}(\xv^{(i)}), \dv = B^{(i)} \xv^{(i)}
   \end{array}
   \right\}.
\end{align*}
Note that we can simply choose $s_i = n$ in (\ref{eq:sqarsity_constrainted}) to model the block-structured smooth nonconvex problem without sparsity constraints. 


\begin{definition}[1-coercive function]\cite{hiriart2004fundamentals}\label{defn:1co}
A function $f : \mathbb{R}^n \to (-\infty,+\infty]$ is called \emph{1-coercive} if
\[
    \frac{f(x)}{\|x\|_2} \;\longrightarrow\; +\infty 
    \qquad \text{as } \|x\|_2\to\infty.
\]
Equivalently, for every $M > 0$ there exists $R > 0$ such that
\[
    \|x\|_2 > R \;\Longrightarrow\; f(x) > M \|x\|_2.
\]
\end{definition}
The notion of 1-coercive functions is commonly used in analysis and optimization~\cite{dubois2025frank}. Some example includes strictly convex quadratic functions and univariate polynomials of even degree with a positive leading coefficient.

The existence of an asymptotically tight Lagrangian dual requires two additional assumptions. Assumption \ref{asp:bounded_variance} requires the nonconvexity of each block in (\ref{prob_smooth_sparsity}) is bounded. Assumption \ref{asp:coercive} implies that $\conv\left(\mathcal{P}^{(i)}_{\textbf{s}}\right)$ is closed and shown in Lemma \ref{lem:guranteen_P}.

\begin{asp}
    \label{asp:bounded_variance}
    There exists some $\beta$ such that $\Xi(\mathcal{P}_{\textbf{s}}^{(i)}) \leq \beta,\forall i \in [k]$.
\end{asp}

A function $f(\xv):\R^n \to (-\infty,+\infty]$ is called proper if there exists some $\xv_0$ such that $f(\xv_0) < \infty$.

\begin{asp}
    \label{asp:coercive}
    Each $f^{(i)}(\cdot)$ is $1$-coercive, closed and proper and $f^{(i)}(\textbf{0}) < \infty$.
\end{asp}

\begin{lemma}
    \label{lem:guanranteen}
    Under Assumption \ref{asp:coercive}, $\conv\left(\mathcal{X}^{(i)}_{\textbf{s}}\right)$ is closed and pointed. Moreover, $\text{rec}\left(\conv\left(\mathcal{X}_\textbf{s}^{(i)}\right)\right)=\{(t,\text{0}): t \geq 0\}$. 
    \begin{proof}
    Let $\chi_{s_i}(\xv) := \begin{cases}
        \infty & \text{ if } \norm{\xv^{(i)}}_0 \geq s_i {+1} \\
        0  & \text{ if otherwise}
    \end{cases}.$ Let $\Tilde{f}^{(i)}(\xv) := f^{(i)}(\xv) + \chi_{s_i}(\xv)$. 
    By the definition of $1$-coercive, it is straightforward that $\Tilde{f}^{(i)}(\xv)$ is $1$-coercive and closed. Since $f^{(i)}(\textbf{0}) < \infty$, $\Tilde{f}^{(i)}(\xv)$ is also proper.
    Moreover, it is clear that $\epi \Tilde{f}^{(i)}(\xv) = \mathcal{X}^{(i)}_{\textbf{s}}$. It is well known that the convex hull of the epigraph of a proper, $1$-coercive and closed function is closed and therefore $\conv\left(\mathcal{X}^{(i)}_{\textbf{s}}\right)$ is closed \cite{hiriart1993convex}. Suppose, for the sake of contradiction, that $\conv\left(\mathcal{X}^{(i)}_{\textbf{s}}\right)$ is not pointed. We first observe that the definition of $1$-coercive-ness implies that there exists some $lb$ such that $lb \leq \Tilde{f}^{(i)}(\xv)$. This implies that $t \geq lb$, for all $(t,\xv) \in \conv\left(\mathcal{X}_\textbf{s}^{(i)}\right)$. Therefore, any line contained in $\conv\left(\mathcal{X}_\textbf{s}^{(i)}\right)$ must take the form of $l(\gamma) := (t_0,\xv_0) + \gamma (0,\rv), \ \gamma \in \R^n$ for some nonzero $\rv$. To also establish the properties of $\text{rec}\left(\conv\left(\mathcal{X}_\textbf{s}^{(i)}\right)\right)$, we show the stronger statement that $\conv\left(\mathcal{X}_\textbf{s}^{(i)}\right)$ does not contain any ray taking the form of $l^+(\gamma) := (t_0,\xv_0) + \gamma (0,\rv), \ \gamma \in \R_+$ for some nonzero $\rv$. This implies that $\conv\left(\mathcal{X}_\textbf{s}^{(i)}\right)$ is pointed.
    Combined with the fact that the function is also lower bounded by some constant, this also implies that $\text{rec}\left(\conv\left(\mathcal{X}_\textbf{s}^{(i)}\right)\right)=\{(t,\text{0}): t \geq 0\}$. The argument is as follows.
    
   Suppose such ray $l^+(\gamma)$ exists, since $\Tilde{f}^{(i)}(\xv)$ is $1$-coercive, there exists some $R > 1$ such that $\norm{\xv}_2 \geq R \implies \Tilde{f}^{(i)}(\xv) \geq (|lb|+|t_0|+1) \norm{\xv}_2$. Choose $|\gamma|$ sufficiently large so that $\norm{\xv_0+ \gamma \rv}_2 \geq (n+2)R$, since $(t_0,\xv_0+ \gamma \rv) \in \conv\left(\mathcal{X}_\textbf{s}^{(i)}\right)$, this implies that $\blue{(t_0,\xv_0+\gamma \rv)}= \sum_{j \in [n+2]}\lambda_j (t_j,\xv_j)$ for some $(t_j,\xv_j) \in \mathcal{X}_{\textbf{s}}^{(i)}$ and $\lambda$ belongs to a simplex of dimension $n+2$. Since $\norm{\cdot}_2$ is subadditive and positively homogeneous, this implies that $\norm{\xv_0+ \gamma \rv }_2 \leq \sum_{i \in [n+2]} \lambda_{j} \norm{\xv_j}_2$.
Since $\norm{\xv_0+ \gamma \rv }_2 \geq (n+2) R$, a simple pigeonhole argument implies that there exists some $j_* \in [n]$, such that
$\lambda_{j_*} \norm{\xv_{j_*}}_2  \geq R$. Since $\lambda_{j_*} \in [0,1]$, this further shows that $\norm{\xv_{j_*}}_2  \geq R$ and by the choice of $R$, this implies that $t_{j_*} \geq (|lb|+|t_0|+1) \norm{\xv_{j_*}}_2$.
This shows that 
\begin{align*}
   \lambda_{j_*} \norm{\xv_{j_*}}_2  \geq R \implies \lambda_{j_*} t_{j_*} \geq \dfrac{R}{\norm{\xv_{j_*}}_2} t_{j_*} \geq R (|lb|+|t_0|+1) \geq |lb|+|t_0|+1.
\end{align*}
This implies that 
\begin{align*}
    \sum_{j \in [n+2]}  \lambda_j t_j =  \lambda_{j_*} t_{j_*} + \sum_{j \in [n+2] \setminus \{{j_*}\}}   \lambda_j t_j  \geq |lb|+|t_0|+1 + (1-\lambda_{j_*}) lb > t_0,
\end{align*}
which leads to \blue{a} contraction.
\end{proof}
\end{lemma}


\begin{lemma}[Theorem 9.1 in \cite{rockafellar2015convex}]
    \label{lem_closed_linear_image}
    Let $C$ be a closed pointed convex set of $\R^n$. For any linear map $\mathcal{A}: \R^n \to \R^m$, if $\text{rec}(C) \cap \ker(\mathcal{A}) = \{0\}$, then $\mathcal{A}(C)$ is closed and pointed. 
\end{lemma}

\begin{lemma}
\label{lem:guranteen_P}
Under assumption \ref{asp:coercive}, $\conv\left(\mathcal{P}^{(i)}_{\textbf{s}}\right)$ is closed and pointed.
\begin{proof}
    By \blue{Lemma} \ref{claim_convex_and_proj}, $\conv\left(\mathcal{P}^{(i)}_{\textbf{s}}\right)$ can be viewed as the image of $\conv\left(\mathcal{X}^{(i)}_{\textbf{s}}\right)$ under the following linear map $\mathcal{A}$: $
         (t,\xv) \to \begin{bmatrix}
       1 &  \\
        & B
    \end{bmatrix} (t,\xv)$.
In this case, we have $\text{rec}\left(\conv\left(\mathcal{X}^{(i)}_{\textbf{s}}\right)\right) \cap \ker(\mathcal{A}) = \{\textbf{0}\}$ because any non-zero $\vv \in \text{rec}\left(\conv\left(\mathcal{X}^{(i)}_{\textbf{s}}\right)\right)$ has form of $(r,\textbf{0})$ for some $r  > 0$ by Lemma \ref{lem:guanranteen} and it is clear that $\mathcal{A} \vv \neq \textbf{0}$ in this case. Applying Lemma \ref{lem_closed_linear_image} and Lemma \ref{lem:guanranteen} yields the result.
\end{proof}
\end{lemma}


The last ingredient for applying Theorem \ref{thm_rem_convex_ball} is to show that there exists a 
smooth hidden convex component ($l_2$-ball) associated with each extreme point of 
$\conv\bigl( \mathcal{P}_{\mathbf{s}}^{(i)} \bigr)$. We obtain this from the following observation. 
The set $\mathcal{P}_{\mathbf{s}}^{(i)}$ can be viewed as a linear projection of a union of epigraphs of certain smooth functions. By the definition of smoothness, for any point in such an epigraph we can find a ball that lies entirely in the epigraph and contains this point, as established in Claim \ref{claim:inscribed_ball} below. Under the linear projection induced by the coupling matrix $B^{(i)}$, this ball is mapped to an ellipsoid in $\mathcal{P}_{\mathbf{s}}^{(i)}$.
This ellipsoid provides a desired smooth hidden convex component. We then extract an $l_2$-ball contained within this ellipsoid, which ultimately serves as the hidden convex component used in our proof.
The maximal distortion under this linear projection and the choice $l_2$-ball contained in this ellipsoid is quantified by $\mathcal{L}_s(B^{(i)})$, defined below. 

\begin{definition}
    Given a matrix $B \in \R^{m \times n}$ and let $\Tilde{B}_{S} \in \R^{(m+1) \times (|S|+1)} := \begin{bmatrix}
        1 & \textbf{0} \\
         \textbf{0} &  B_{S} 
    \end{bmatrix}$ be a matrix 
    where $B_S$ is the submatrix of $B$ that consists of columns of $B$ with indices in $S$.
    We define the \textit{projection factor with respect to sparsity level $s$} as
    \begin{align*}
    \mathcal{L}_s(B) := \inf_{S \subseteq [n]:|S| = s} \sigma_{\diamond}\left( \Tilde{B}_S (\Tilde{B}_S)^{\top}   \right)
    \end{align*}
where $\sigma_{\diamond}(A) := \begin{cases}
    0 & \text{ if } A \text{ is not positive definite}; \\ 
    \frac{{\sigma_{\min}(A^{-1/2})}}{\sigma_{\max}^2(A^{-1/2}) } & \text{ otherwise};
\end{cases}$ and $\sigma_{\min}(\cdot),\sigma_{\max}(\cdot)$ are the minimal singular value and maximal singular value, respectively. $A^{1/2}$ denotes the principal square root of a positive definite matrix $A$.
\end{definition}

\begin{remark}
    To obtain a nondegenerate hidden convex component, 
we require $\mathcal{L}_s(B^{(i)}) > 0$. This condition is equivalent to the existence of 
$m$ linear independent columns among every possible collection of $s^{(i)}$ columns of $B^{(i)}$.
A necessary condition for this is $s^{(i)} \ge m$. Apart from this simple dimensional 
requirement, the condition $\mathcal{L}_s(B^{(i)}) > 0$ is not very restrictive: it holds 
generically and can be ensured almost surely by adding an arbitrarily small smooth random 
perturbation to $B^{(i)}$ whenever $s^{(i)} \ge m$.
\end{remark}


\blue{
\begin{definition}
A differentiable function $f:\R^n \to \R$ is called locally $L$-smooth with respect to a set $M \subseteq \R^n$ if
$f(\yv) \leq f(\xv) + \iprod{\yv-\xv}{\nabla f(\xv)} + \dfrac{L}{2} \norm{\yv-\xv}_2^2,\forall \yv,\xv \in M$.
\end{definition}
}

\blue{
\begin{claim}
    \label{claim:inscribed_ball}
    Let $M\subseteq\mathbb{R}^n$ and let $f:\mathbb{R}^n\to\mathbb{R}$ be $L$-smooth on $M+\mathcal{B}_2(\mathbf{0},2/L)$. Then for every $\xv\in M$ there is an $\ell_2$-ball $\mathcal{B}\subseteq\mathbb{R}^{n+1}$ of radius $1/L$ with $ (\xv,f(\xv))\in\mathcal{B}\subseteq\operatorname{epi} f.$
\end{claim}
}
\begin{proof}
Fix a \(\xv\in M\). Let $\gv:=\nabla f(\xv) \text{ and }
    s:=\sqrt{1+\|\gv\|_2^2}.$
Define
\[
    \cv_{\xv}
    :=
    (\xv,f(\xv))
    +
    \frac{1}{L}\frac{(-\gv,1)}{\sqrt{1+\|\gv\|_2^2}}
    =
    \left(
        \xv-\frac{\gv}{Ls},
        f(\xv)+\frac{1}{Ls}
    \right).
\]
We claim that
$
    \mathcal{B}_{\xv}:=\mathcal{B}_2\left(\cv_{\xv},\frac{1}{L}\right)
    \subseteq \operatorname{epi} f.
$
First, observe that
$
    \left\|(\xv,f(\xv))-\cv_{\xv}\right\|_2
    =
    \left\|
        \left(
            \frac{\gv}{Ls},
            -\frac{1}{Ls}
        \right)
    \right\|_2  
    =
    \frac{1}{Ls}\sqrt{\|\gv\|_2^2+1}
    =
    \frac{1}{L}.
$
Therefore $(\xv,f(\xv))\in \mathcal{B}_{\xv}.$

Now take any \((\yv,t)\in \mathcal{B}_{\xv}\), and set $  \hv:=\yv-\xv.$
Since \((\yv,t)\in \mathcal{B}_{\xv}\), we have $\left\|
        \left(
            \hv+\frac{\gv}{Ls},
            t-f(\xv)-\frac{1}{Ls}
        \right)
    \right\|_2
    \le \frac{1}{L}.$
This implies $ \left\|\hv+\frac{\gv}{Ls}\right\|_2\le \frac{1}{L}.$ This also implies that $\|\hv\|_2
    \le
    \left\|\hv+\frac{\gv}{Ls}\right\|_2
    +
    \frac{\|\gv\|_2}{Ls}  
    \le
    \frac{1}{L}
    +
    \frac{1}{L}
    =
    \frac{2}{L}$
since \(\|\gv\|_2/s\le 1\). Since \(\xv\in M\), this implies $\yv=\xv+\hv\in M+\mathcal{B}_2(\mathbf{0},2/L).$
Therefore, by local \(L\)-smoothness of \(f\) on \(M+\mathcal{B}_2(\mathbf{0},2/L)\), we have $f(\yv)
    \le
    f(\xv)+\langle \gv,\hv\rangle+\frac{L}{2}\|\hv\|_2^2.$

It remains to show that the right-hand side is at most \(t\) and this will imply that $(\yv,t) \in \operatorname{epi} f$. Since
\((\yv,t)\in \mathcal{B}_{\xv}\), we have $ \left(
        t-f(\xv)-\frac{1}{Ls}
    \right)^2
    +
    \left\|\hv+\frac{\gv}{Ls}\right\|_2^2
    \le
    \frac{1}{L^2}.$
Thus $ t
    \ge
    f(\xv)+\frac{1}{Ls}
    -
    \sqrt{
        \frac{1}{L^2}
        -
        \left\|\hv+\frac{\gv}{Ls}\right\|_2^2
    }.$
The square root is well-defined because $ \left\|\hv+\frac{\gv}{Ls}\right\|_2\le \frac{1}{L},$
and hence $    \frac{1}{L^2}
    -
    \left\|\hv+\frac{\gv}{Ls}\right\|_2^2
    \ge 0.$
We now prove the elementary inequality
\[
    \langle \gv,\hv\rangle+\frac{L}{2}\|\hv\|_2^2
    \le
    \frac{1}{Ls}
    -
    \sqrt{
        \frac{1}{L^2}
        -
        \left\|\hv+\frac{\gv}{Ls}\right\|_2^2
    }.
\]
Equivalently, it is enough to show
\[
    \sqrt{
        \frac{1}{L^2}
        -
        \left\|\hv+\frac{\gv}{Ls}\right\|_2^2
    }
    \le
    \frac{1}{Ls}
    -
    \langle \gv,\hv\rangle
    -
    \frac{L}{2}\|\hv\|_2^2.
\]

First, the right-hand side is nonnegative. Indeed, from $\left\|\hv+\frac{\gv}{Ls}\right\|_2\le \frac{1}{L},$
we get $\|\hv\|_2^2
    +
    \frac{2}{Ls}\langle \gv,\hv\rangle
    +
    \frac{\|\gv\|_2^2}{L^2s^2}
    \le
    \frac{1}{L^2}.$
Since \(s^2=1+\|\gv\|_2^2\), this gives $ \|\hv\|_2^2
    +
    \frac{2}{Ls}\langle \gv,\hv\rangle
    \le
    \frac{1}{L^2s^2}.$
Multiplying both sides by \(Ls/2>0\), we obtain $\langle \gv,\hv\rangle
    +
    \frac{Ls}{2}\|\hv\|_2^2
    \le
    \frac{1}{2Ls}.$
Therefore
\[
\begin{aligned}
    \frac{1}{Ls}
    -
    \langle \gv,\hv\rangle
    -
    \frac{L}{2}\|\hv\|_2^2
    &=
    \frac{1}{Ls}
    -
    \left(
        \langle \gv,\hv\rangle
        +
        \frac{Ls}{2}\|\hv\|_2^2
    \right)
    +
    \frac{L(s-1)}{2}\|\hv\|_2^2  \\
    &\ge
    \frac{1}{Ls}
    -
    \frac{1}{2Ls}
    +
    \frac{L(s-1)}{2}\|\hv\|_2^2  \\
    &\ge 0.
\end{aligned}
\]

Now both sides are nonnegative, so it suffices to square both sides. We compute
\[
\begin{aligned}
&\left(
    \frac{1}{Ls}
    -
    \langle \gv,\hv\rangle
    -
    \frac{L}{2}\|\hv\|_2^2
\right)^2
-
\left(
    \frac{1}{L^2}
    -
    \left\|\hv+\frac{\gv}{Ls}\right\|_2^2
\right)  \\
&=
\left(
    \langle \gv,\hv\rangle
    +
    \frac{L}{2}\|\hv\|_2^2
\right)^2
+
\left(1-\frac{1}{s}\right)\|\hv\|_2^2
\ge 0.
\end{aligned}
\]
Hence
\[
    \sqrt{
        \frac{1}{L^2}
        -
        \left\|\hv+\frac{\gv}{Ls}\right\|_2^2
    }
    \le
    \frac{1}{Ls}
    -
    \langle \gv,\hv\rangle
    -
    \frac{L}{2}\|\hv\|_2^2.
\]
Thus
\[
    \langle \gv,\hv\rangle+\frac{L}{2}\|\hv\|_2^2
    \le
    \frac{1}{Ls}
    -
    \sqrt{
        \frac{1}{L^2}
        -
        \left\|\hv+\frac{\gv}{Ls}\right\|_2^2
    }.
\]
Combining this with the lower bound on \(t\), we get
\[
\begin{aligned}
    t
    &\ge
    f(\xv)+\frac{1}{Ls}
    -
    \sqrt{
        \frac{1}{L^2}
        -
        \left\|\hv+\frac{\gv}{Ls}\right\|_2^2
    } \\
    &\ge
    f(\xv)+\langle \gv,\hv\rangle+\frac{L}{2}\|\hv\|_2^2 \\
    &\ge
    f(\yv).
\end{aligned}
\]
Therefore \((\yv,t)\in \operatorname{epi} f\). Since \((\yv,t)\in \mathcal{B}_{\xv}\) was arbitrary, we conclude that $ \mathcal{B}_{\xv}\subseteq \operatorname{epi} f.$
\end{proof}

\begin{lemma}
\label{lem_liu_gauge}
\blue{(Table 1 and Proposition 3 in \cite{liu2023gauges})}
    Let $E:=\{\xv: \norm{Ax-\kv}_2 \leq 1\}$ be an ellipsoid, where $A$ is an invertible matrix. Then, for every $\pv \in E$, there exists a $l_2$-norm ball containing $\pv$ that is subset of $E$ and has a radius of 
     $\frac{\sigma_{\min}(A)}{\sigma_{\max}(A)^2}.$
\end{lemma}

\begin{claim} 
    \label{claim_existance_of_ball_proj}
    \blue{Let $\vv := \left(f(\xv_{\vv}), B^{(i)} \xv_{\vv}\right)$ be an extreme point of $\conv\left(\mathcal{P}_\textbf{s}^{(i)}\right)$, where $\xv_{\vv} \in \R^n$. If $f$ is $L$-smooth on $\mathcal{B}_2\left(0, \norm{\xv_{\vv}}_2 + \frac{2}{L}\right)$, then there exists an $\ell_2$-ball $\mathcal{B}$ of radius $\frac{\mathcal{L}_s(B^{(i)})}{L}$ such that $\vv \in \mathcal{B} \subseteq \mathcal{P}_\textbf{s}^{(i)}$.}
    
\begin{proof}
For a fixed support set $S \subseteq [n]$, we define
\begin{equation}
\begin{aligned}
    \mathcal{P}^{(i)}_{S} := \left\{ (t,\dv) \in \R^{m+1} \;\middle\vert\;
   \begin{array}{@{}l@{}} \exists \xv, x_i = 0,\forall i \notin S; \\ \dv = B^{(i)} \xv;\\
   t \geq f(\xv);
   \end{array} 
\right\},
\end{aligned}
\end{equation}
It is clear that
\begin{align*}
    \mathcal{P}^{(i)}_\textbf{s} = \bigcup_{S \subseteq [n], |S| = s_i } \mathcal{P}^{(i)}_S.
\end{align*}

Let $\vv:=(t_{\vv},\dv_{\vv})$ be any extreme point of $\conv( \mathcal{P}_\textbf{s}^{(i)})$. In this case, there must exist some 
$S \subseteq [n]$, with $|S| = s_i$, such that $\vv \in P^{(i)}_S$. By construction of $ P^{(i)}_S$, there exists some $\xv_{\vv} \in \R^n$ such that 
\begin{align*}
    & \dv_{\vv} = B^{(i)}\xv_{\vv}, t_{\vv} \geq f^{(i)}(\xv_{\vv}), \\
    & (\xv_{\vv})_i =0,\forall i \notin S.
\end{align*}
Since $\vv$ is a extreme point, we can further assume that $t_{\vv} = f^{(i)}(\xv_\vv)$. Consider the following function $f^{(i)}_S(\yv):\R^S \to \R$ which is constructed from $f^{(i)}(\xv)$ by fixing the variables outside of $S$ to be zero. \blue{Clearly, $f^{(i)}_S$ is a locally $L$-smooth convex function with respect to $\mathcal{B}_2\left(0,\norm{(\xv_{\vv})_{S}}_2+\frac{2}{L}\right)$ where $(\xv_{\vv})_{S}$ is the subvector of $\xv_{\vv}$ indexed by $S$}. Applying Claim \ref{claim:inscribed_ball}, there is a ball $\mathcal{B}_{S}$ with radius $\frac{1}{L}$ with dimension $|S|+1$ that $(f(\xv_{\vv}), (\xv_{\vv})_S)\subseteq \mathcal{B}_{S}\subseteq \epi f_S^{(i)}$. Now consider the linear map:
\begin{align*}
    \blue{(t,\xv_S) \to \begin{bmatrix}
       1 & \textbf{0} \\
       \textbf{0} & B_S^{(i)}
    \end{bmatrix} (t,\xv_S)}
\end{align*}
\blue{Let $Q :=
    \begin{bmatrix}
       1 & \mathbf{0} \\
       \mathbf{0} & B_S^{(i)}
    \end{bmatrix}.$
Then the above linear map can be written compactly as $(t,\xv_S)\mapsto Q(t,\xv_S).$
If \(QQ^\top\) is not invertible, then by the definition of $\mathcal{L}_s(B^{(i)})$, we have $\mathcal{L}_s(B^{(i)})=0$.
In this degenerate case, the desired ball has radius \(0\), and the claim follows directly since $ \{\vv\}\subseteq \mathcal{P}^{(i)}_S\subseteq \mathcal{P}^{(i)}_{\mathbf{s}}.$
Therefore, we may assume that \(QQ^\top\) is invertible.
Let \(\bar{\zv}\) denote the center of the ball \(\mathcal{B}_S\). Since
\(\mathcal{B}_S\) has radius \(1/L\), we can write $\mathcal{B}_S
    =
    \left\{
        \zv :
        \|\zv-\bar{\zv}\|_2\le \frac{1}{L}
    \right\}.$
The image of \(\mathcal{B}_S\) under the linear map \(Q\) is $ E:=Q\mathcal{B}_S.$
Because \(QQ^\top\) is invertible, this image is an ellipsoid in
\(\mathbb{R}^{m+1}\). More explicitly, $E
    =
    \left\{
        \yv :
        (\yv-Q\bar{\zv})^\top
        (QQ^\top)^{-1}
        (\yv-Q\bar{\zv})
        \le
        \frac{1}{L^2}
    \right\}.$
Equivalently, if $ U:=(QQ^\top)^{-1/2}$, then $E
    =
    \left\{
        \yv :
        \|U(\yv-Q\bar{\zv})\|_2
        \le
        \frac{1}{L}
    \right\}.$
Since $\mathcal{B}_S\subseteq \operatorname{epi} f_S^{(i)}$,
we have $E=Q\mathcal{B}_S
    \subseteq
    Q\operatorname{epi} f_S^{(i)}
    =
    \mathcal{P}^{(i)}_S
    \subseteq
    \mathcal{P}^{(i)}_{\mathbf{s}}.$
Moreover, as shown above, $\vv
    =
    Q\big(f(\xv_{\vv}),(\xv_{\vv})_S\big)
    \in E.$
By Lemma \ref{lem_liu_gauge}, together with the definition of
\(\mathcal{L}_s(B^{(i)})\), the ellipsoid \(E\) contains an \(\ell_2\)-ball
\(\mathcal{B}\) of radius $ \frac{\mathcal{L}_s(B^{(i)})}{L}$.
}

\end{proof}

\end{claim}




\blue{
As mentioned earlier, the Lagrangian dual of \eqref{prob_smooth_sparsity} admits the following equivalent primal characterization:
\begin{equation}
    \label{eq:dual_sqarsity_constrainted}
    \begin{aligned}
         \DUAL_{\textbf{s}}(\bv) := \inf & \sum_{i=1}^{k} t^{(i)} \\
    \text{s.t. } & \sum_{i=1}^k B^{(i)} \xv^{(i)} \leq \bv, \\
    & (t^{(i)},\xv^{(i)}) \in \conv\left(\mathcal{X}^{(i)}_{\textbf{s}}\right).  \\
    \end{aligned}
\end{equation}
}
\blue{
For simplicity, we assume that (\ref{eq:dual_sqarsity_constrainted}) attains its optimal solution (otherwise one can replace the argument below with an $\epsilon$-optimal solution).
Using Lemma \label{lem:guranteen_P} and Lemma \ref{lem:face_decompose}, we can find an optimal solution $\left\{ (\tv^{(i)}_{*},\xv^{(i)}_{*})\right\}_{i=1}^k$ of (\ref{eq:dual_sqarsity_constrainted}) and faces $\{F^{(i)}\}_{i=1}^k$ of $\conv\left(\mathcal{P}^{(i)}_{\textbf{s}}\right)$ such that (i) $\sum_{i=1}^k (t_*^{(i)}, B^{(i)} \xv^{(i)}_*) \in \sum F^{(i)}$ and (ii) there exists an index set $\mathcal{I} := \{i \in[k]:\dim(F^{(i)})>0\}$ such that $|\mathcal{I}|\leq m+1$.
}

\begin{theorem}
\label{thm_application}
  \blue{  Assume Assumption \ref{asp:strong_dual}, Assumption \ref{asp:bounded_variance} and Assumption \ref{asp:coercive}.
    Consider the problem (\ref{eq:dual_sqarsity_constrainted}), let $\left\{ (\tv^{(i)}_{*},\xv^{(i)}_{*})\right\}_{i=1}^k$ be the optimal solution discussed above. Suppose for each $i \in [k] \setminus \mathcal{I}$, $f^{(i)}(\cdot)$ is locally $L^{(i)}$-smooth with respect to the set $ \mathcal{B}_2\left(\textbf{0},\norm{\xv_*^{(i)}}_2+\frac{2}{L^{(i)}}\right)$. Then
        \begin{align*}
         & \OPT_{\textbf{s}}(\bv+\mathcal{E} \textbf{1}) - \mathcal{E} \leq \DUAL_{\textbf{s}}(\bv) \leq \OPT_{\textbf{s}}(\bv),
    \end{align*}
    where
    \begin{align*}
                & \mathcal{E} \leq \sqrt{\left( \sum_{i \in [k] \setminus I} \frac{ \mathcal{L}_s(B^{(i)}) }{L^{(i)}} \right)^2 +   (m+1) \beta^2 }- \left( \sum_{i \in [k] \setminus I} \frac{ \mathcal{L}_s(B^{(i)}) }{L^{(i)}} \right). 
    \end{align*}
If there exist two universal constants $\omega_1,\omega_2$ such that $\mathcal{L}_{s}(B^{(i)}) \geq \omega_1, L^{(i)} \leq \omega_2, \forall i \in [k]$, it then follows that
   \begin{align*}
        \mathcal{E}  \to 0 \text{ as } k \to \infty.
   \end{align*}
   }
\end{theorem}

    
   \begin{proof}
Using Claim \ref{claim_existance_of_ball_proj}, we can equip each $(t_*^{(i)}, B^{(i)} \xv^{(i)}_*)$ for $i \in [k] \setminus \mathcal{I}$ with a $l_2$-ball with radius $\frac{ \mathcal{L}_s(B^{(i)}) }{L}$ as a hidden convex component. The remaining statement is then direct consequence of Lemma \ref{lem:local_geom} and Theorem \ref{thm_duality_gap}.
\end{proof}

\begin{example}
\blue{
    Consider an instance of (\ref{eq:sqarsity_constrainted})  where for all $i \in [k]$, we have $f^{(i)} = \sum_{j \in [n]} (x_j-1)^4-(x_j-1)^2$, $\norm{B^{(i)}}_{\text{op}} \leq 1$,  $\mathcal{L}_{s}(B^{(i)}) \geq \omega_1$ with some constant $\omega_1$, $s^{(i)}=n$ (no sparsity constraint) and
    $\bv \geq 0$.
    Each $f^{(i)}$ is 1-coercive and $\Xi(\epi f) \leq \sqrt{n/2}$. In this case, we have $\Xi(\mathcal{P}_{s}^{(i)})\leq \norm{B^{(i)}}_{\text{op}} \Xi(\epi f)  \leq \sqrt{n/2}$. Let $(\xv_*^{(i)},t_*^{(i)})$ be the optimal solution of dual problem (\ref{eq:dual_sqarsity_constrainted}). We first observe that at least half of optimal solution must be bounded due to 1-coercivity. Indeed, we have $(\textbf{0},0)$ a feasible solution. Therefore, we must have
$\sum_{i=1}^k t_*^{(i)} \leq 0.$
Therefore, the average objective value of $t_*^{(i)}$ is at least $0$. Since each $t_*^{(i)}$ is also lower bounded by $-0.25n$, applying the Markov inequality or a basic counting argument, there exists an index set $O$ with $|O| \geq 0.5 k$ such that $i \in O \implies t^{(i)}_* \leq 0.25 n$. Since $f^{(i)}$ is 1-coercive, there exists a constant $r_1$ (independent of $k$) such that $i \in O \implies \norm{\xv_*^{(i)}}_2 \leq r_1$. In this case, we can find a constant $0 < r_2 < \infty$ (independent of $k$) such that $i \in O \implies \text{$f^{(i)}(\cdot)$ is locally $r_2$-smooth with respect to the set $ \mathcal{B}_2\left(\textbf{0},\norm{\xv_*^{(i)}}_2+\frac{2}{r_2}\right)$}$. For those $i \in [k] \setminus O$, we simply set their smoothness constant $L^{(i)} = \infty$. Applying Theorem \ref{thm_application}, we have
\begin{align*}
    \OPT_{\textbf{s}}(\bv+\mathcal{E} \textbf{1}) - \mathcal{E} \leq \DUAL_{\textbf{s}}(\bv) \leq \OPT_{\textbf{s}}(\bv)
\end{align*}
where
    \begin{align*}
                & \mathcal{E} \leq \sqrt{\left( \frac{w_1}{r_2} \left(\frac{k}{2}-m-1\right) \right)^2 +   (m+1) \frac{n}{2} }- \left( \frac{w_1}{r_2} \left(\frac{k}{2}-m-1\right) \right). 
    \end{align*}
Clearly, $\mathcal{E}  \to 0 \text{ as } k \to \infty$.
}
\end{example}

\begin{remark}
   \blue{ The drawback of Theorem \ref{thm:generalrefinement} and Theorem \ref{thm_application} is that the required smoothness structure is generally not preserved under intersection. For example, although each \(\ell_2\)-ball has a smooth boundary, the intersection of two \(\ell_2\)-balls may have a nonsmooth boundary at the points where the two spherical boundaries meet. Therefore, the current framework does not directly allow constraints of the form \(\xv^{(i)}\in \mathcal{X}^{(i)}\) for a general feasible set \(\mathcal{X}^{(i)}\) in \((\ref{eq:sqarsity_constrainted})\).
One possible way to address this issue is to replace the constrained problem by a smooth approximation:
\[
    \min_{\xv\in \mathcal{X}^{(i)}} f(\xv)
    \quad \longrightarrow \quad
    \min_{\xv} \; f(\xv)+g_{\mathcal{X}^{(i)}}(\xv),
\]
where \(g_{\mathcal{X}^{(i)}}\) is a smooth barrier or penalty function associated with \(\mathcal{X}^{(i)}\). The resulting objective is typically locally smooth when \(f\) is smooth. However, additional assumptions are needed to control the nonconvexity of \(f+g_{\mathcal{X}^{(i)}}\), and the approximation error introduced by the barrier or penalty term must also be taken into account.}
\end{remark}



\section*{Acknowledgments}
We would also like to thank Diego Cifuentes, Greg Blekherman, Santosh Vempala, and Kazuo Murota for the various wonderful discussions.



\appendix


\bibliography{ref}


\begin{thebibliography}{27}
\ifx \bisbn   \undefined \def \bisbn  #1{ISBN #1}\fi
\ifx \binits  \undefined \def \binits#1{#1}\fi
\ifx \bauthor  \undefined \def \bauthor#1{#1}\fi
\ifx \batitle  \undefined \def \batitle#1{#1}\fi
\ifx \bjtitle  \undefined \def \bjtitle#1{#1}\fi
\ifx \bvolume  \undefined \def \bvolume#1{\textbf{#1}}\fi
\ifx \byear  \undefined \def \byear#1{#1}\fi
\ifx \bissue  \undefined \def \bissue#1{#1}\fi
\ifx \bfpage  \undefined \def \bfpage#1{#1}\fi
\ifx \blpage  \undefined \def \blpage #1{#1}\fi
\ifx \burl  \undefined \def \burl#1{\textsf{#1}}\fi
\ifx \doiurl  \undefined \def \doiurl#1{\url{https://doi.org/#1}}\fi
\ifx \betal  \undefined \def \betal{\textit{et al.}}\fi
\ifx \binstitute  \undefined \def \binstitute#1{#1}\fi
\ifx \binstitutionaled  \undefined \def \binstitutionaled#1{#1}\fi
\ifx \bctitle  \undefined \def \bctitle#1{#1}\fi
\ifx \beditor  \undefined \def \beditor#1{#1}\fi
\ifx \bpublisher  \undefined \def \bpublisher#1{#1}\fi
\ifx \bbtitle  \undefined \def \bbtitle#1{#1}\fi
\ifx \bedition  \undefined \def \bedition#1{#1}\fi
\ifx \bseriesno  \undefined \def \bseriesno#1{#1}\fi
\ifx \blocation  \undefined \def \blocation#1{#1}\fi
\ifx \bsertitle  \undefined \def \bsertitle#1{#1}\fi
\ifx \bsnm \undefined \def \bsnm#1{#1}\fi
\ifx \bsuffix \undefined \def \bsuffix#1{#1}\fi
\ifx \bparticle \undefined \def \bparticle#1{#1}\fi
\ifx \barticle \undefined \def \barticle#1{#1}\fi
\bibcommenthead
\ifx \bconfdate \undefined \def \bconfdate #1{#1}\fi
\ifx \botherref \undefined \def \botherref #1{#1}\fi
\ifx \url \undefined \def \url#1{\textsf{#1}}\fi
\ifx \bchapter \undefined \def \bchapter#1{#1}\fi
\ifx \bbook \undefined \def \bbook#1{#1}\fi
\ifx \bcomment \undefined \def \bcomment#1{#1}\fi
\ifx \oauthor \undefined \def \oauthor#1{#1}\fi
\ifx \citeauthoryear \undefined \def \citeauthoryear#1{#1}\fi
\ifx \endbibitem  \undefined \def \endbibitem {}\fi
\ifx \bconflocation  \undefined \def \bconflocation#1{#1}\fi
\ifx \arxivurl  \undefined \def \arxivurl#1{\textsf{#1}}\fi
\csname PreBibitemsHook\endcsname

\bibitem[\protect\citeauthoryear{L{\"u}bbecke}{2010}]{lubbecke2010column}
\begin{barticle}
\bauthor{\bsnm{L{\"u}bbecke}, \binits{M.E.}}:
\batitle{Column generation}.
\bjtitle{Wiley encyclopedia of operations research and management science}
\bvolume{17},
\bfpage{18}--\blpage{19}
(\byear{2010})
\end{barticle}
\endbibitem

\bibitem[\protect\citeauthoryear{Desaulniers et~al.}{2006}]{desaulniers2006column}
\begin{bbook}
\bauthor{\bsnm{Desaulniers}, \binits{G.}},
\bauthor{\bsnm{Desrosiers}, \binits{J.}},
\bauthor{\bsnm{Solomon}, \binits{M.M.}}:
\bbtitle{Column Generation}
vol. \bseriesno{5}.
\bpublisher{Springer},
\blocation{New York, USA}
(\byear{2006})
\end{bbook}
\endbibitem

\bibitem[\protect\citeauthoryear{Cifuentes et~al.}{2025}]{cifuentes2025lagrangian}
\begin{bchapter}
\bauthor{\bsnm{Cifuentes}, \binits{D.}},
\bauthor{\bsnm{Dey}, \binits{S.S.}},
\bauthor{\bsnm{Xu}, \binits{J.}}:
\bctitle{Lagrangian dual for integer optimization with zero duality gap that admits decomposition}.
In: \bbtitle{International Conference on Integer Programming and Combinatorial Optimization},
pp. \bfpage{184}--\blpage{198}
(\byear{2025}).
\bcomment{Springer}
\end{bchapter}
\endbibitem

\bibitem[\protect\citeauthoryear{Starr}{1969}]{starr1969quasi}
\begin{botherref}
\oauthor{\bsnm{Starr}, \binits{R.M.}}:
Quasi-equilibria in markets with non-convex preferences.
Econometrica: journal of the Econometric Society,
25--38
(1969)
\end{botherref}
\endbibitem

\bibitem[\protect\citeauthoryear{Ben-Tal and Nemirovski}{2001}]{ben2001lectures}
\begin{bbook}
\bauthor{\bsnm{Ben-Tal}, \binits{A.}},
\bauthor{\bsnm{Nemirovski}, \binits{A.}}:
\bbtitle{Lectures on Modern Convex Optimization: Analysis, Algorithms, and Engineering Applications}.
\bpublisher{SIAM},
\blocation{Philadelphia, USA}
(\byear{2001})
\end{bbook}
\endbibitem

\bibitem[\protect\citeauthoryear{P{\'o}lik and Terlaky}{2007}]{polik2007survey}
\begin{barticle}
\bauthor{\bsnm{P{\'o}lik}, \binits{I.}},
\bauthor{\bsnm{Terlaky}, \binits{T.}}:
\batitle{A survey of the {S}-lemma}.
\bjtitle{SIAM review}
\bvolume{49}(\bissue{3}),
\bfpage{371}--\blpage{418}
(\byear{2007})
\end{barticle}
\endbibitem

\bibitem[\protect\citeauthoryear{Chandrasekaran et~al.}{2025}]{chandrasekaran2025lagrangian}
\begin{botherref}
\oauthor{\bsnm{Chandrasekaran}, \binits{V.}},
\oauthor{\bsnm{Duff}, \binits{T.}},
\oauthor{\bsnm{Rodriguez}, \binits{J.I.}},
\oauthor{\bsnm{Shu}, \binits{K.}}:
Lagrangian dual sections: A topological perspective on hidden convexity.
arXiv preprint arXiv:2510.06112
(2025)
\end{botherref}
\endbibitem

\bibitem[\protect\citeauthoryear{Dines}{1941}]{dines1941mapping}
\begin{barticle}
\bauthor{\bsnm{Dines}, \binits{L.L.}}:
\batitle{On the mapping of quadratic forms}.
\bjtitle{Bulletin of the American Mathematical Society}
\bvolume{47}(\bissue{6}),
\bfpage{494}--\blpage{498}
(\byear{1941})
\end{barticle}
\endbibitem

\bibitem[\protect\citeauthoryear{Fradelizi et~al.}{2018}]{fradelizi2018convexification}
\begin{barticle}
\bauthor{\bsnm{Fradelizi}, \binits{M.}},
\bauthor{\bsnm{Madiman}, \binits{M.}},
\bauthor{\bsnm{Marsiglietti}, \binits{A.}},
\bauthor{\bsnm{Zvavitch}, \binits{A.}}:
\batitle{The convexification effect of minkowski summation}.
\bjtitle{EMS Surveys in Mathematical Sciences}
\bvolume{5}(\bissue{1}),
\bfpage{1}--\blpage{64}
(\byear{2018})
\end{barticle}
\endbibitem

\bibitem[\protect\citeauthoryear{Bertsekas et~al.}{1983}]{bertsekas1983optimal}
\begin{barticle}
\bauthor{\bsnm{Bertsekas}, \binits{D.}},
\bauthor{\bsnm{Lauer}, \binits{G.}},
\bauthor{\bsnm{Sandell}, \binits{N.}},
\bauthor{\bsnm{Posbergh}, \binits{T.}}:
\batitle{Optimal short-term scheduling of large-scale power systems}.
\bjtitle{IEEE Transactions on Automatic Control}
\bvolume{28}(\bissue{1}),
\bfpage{1}--\blpage{11}
(\byear{1983})
\end{barticle}
\endbibitem

\bibitem[\protect\citeauthoryear{Aubin and Ekeland}{1976}]{aubin1976estimates}
\begin{barticle}
\bauthor{\bsnm{Aubin}, \binits{J.-P.}},
\bauthor{\bsnm{Ekeland}, \binits{I.}}:
\batitle{Estimates of the duality gap in nonconvex optimization}.
\bjtitle{Mathematics of Operations Research}
\bvolume{1}(\bissue{3}),
\bfpage{225}--\blpage{245}
(\byear{1976})
\end{barticle}
\endbibitem

\bibitem[\protect\citeauthoryear{Udell and Boyd}{2016}]{udell2016bounding}
\begin{barticle}
\bauthor{\bsnm{Udell}, \binits{M.}},
\bauthor{\bsnm{Boyd}, \binits{S.}}:
\batitle{Bounding duality gap for separable problems with linear constraints}.
\bjtitle{Computational Optimization and Applications}
\bvolume{64}(\bissue{2}),
\bfpage{355}--\blpage{378}
(\byear{2016})
\end{barticle}
\endbibitem

\bibitem[\protect\citeauthoryear{Bi and Tang}{2020}]{bi2016refined}
\begin{barticle}
\bauthor{\bsnm{Bi}, \binits{Y.}},
\bauthor{\bsnm{Tang}, \binits{A.}}:
\batitle{Duality gap estimation via a refined shapley--folkman lemma}.
\bjtitle{SIAM Journal on Optimization}
\bvolume{30}(\bissue{2}),
\bfpage{1094}--\blpage{1118}
(\byear{2020})
\end{barticle}
\endbibitem

\bibitem[\protect\citeauthoryear{Kerdreux et~al.}{2023}]{kerdreux2023stable}
\begin{barticle}
\bauthor{\bsnm{Kerdreux}, \binits{T.}},
\bauthor{\bsnm{Colin}, \binits{I.}},
\bauthor{\bsnm{d’Aspremont}, \binits{A.}}:
\batitle{Stable bounds on the duality gap of separable nonconvex optimization problems}.
\bjtitle{Mathematics of Operations Research}
\bvolume{48}(\bissue{2}),
\bfpage{1044}--\blpage{1065}
(\byear{2023})
\end{barticle}
\endbibitem

\bibitem[\protect\citeauthoryear{Nguyen and Vohra}{2024}]{nguyen2022near}
\begin{barticle}
\bauthor{\bsnm{Nguyen}, \binits{T.}},
\bauthor{\bsnm{Vohra}, \binits{R.}}:
\batitle{(near-) substitute preferences and equilibria with indivisibilities}.
\bjtitle{Journal of Political Economy}
\bvolume{132}(\bissue{12}),
\bfpage{4122}--\blpage{4154}
(\byear{2024})
\end{barticle}
\endbibitem

\bibitem[\protect\citeauthoryear{Murota and Tamura}{2025}]{murota2025shapley}
\begin{barticle}
\bauthor{\bsnm{Murota}, \binits{K.}},
\bauthor{\bsnm{Tamura}, \binits{A.}}:
\batitle{Shapley--{F}olkman-type theorem for integrally convex sets}.
\bjtitle{Discrete Applied Mathematics}
\bvolume{360},
\bfpage{42}--\blpage{50}
(\byear{2025})
\end{barticle}
\endbibitem

\bibitem[\protect\citeauthoryear{Wegmann}{1980}]{wegmann1980einige}
\begin{barticle}
\bauthor{\bsnm{Wegmann}, \binits{R.}}:
\batitle{Einige ma{\ss}zahlen f{\"u}r nichtkonvexe mengen}.
\bjtitle{Archiv der Mathematik}
\bvolume{34}(\bissue{1}),
\bfpage{69}--\blpage{74}
(\byear{1980})
\end{barticle}
\endbibitem

\bibitem[\protect\citeauthoryear{Hiriart-Urruty and Lemar{\'e}chal}{1993}]{hiriart1993convex}
\begin{botherref}
\oauthor{\bsnm{Hiriart-Urruty}, \binits{J.-B.}},
\oauthor{\bsnm{Lemar{\'e}chal}, \binits{C.}}:
Convex analysis and minimization algorithms ii.
Grundlehren der mathematischen Wissenschaften
(1993)
\end{botherref}
\endbibitem

\bibitem[\protect\citeauthoryear{Boland et~al.}{2018}]{boland2018combining}
\begin{barticle}
\bauthor{\bsnm{Boland}, \binits{N.}},
\bauthor{\bsnm{Christiansen}, \binits{J.}},
\bauthor{\bsnm{Dandurand}, \binits{B.}},
\bauthor{\bsnm{Eberhard}, \binits{A.}},
\bauthor{\bsnm{Linderoth}, \binits{J.}},
\bauthor{\bsnm{Luedtke}, \binits{J.}},
\bauthor{\bsnm{Oliveira}, \binits{F.}}:
\batitle{Combining progressive hedging with a frank--wolfe method to compute lagrangian dual bounds in stochastic mixed-integer programming}.
\bjtitle{SIAM Journal on Optimization}
\bvolume{28}(\bissue{2}),
\bfpage{1312}--\blpage{1336}
(\byear{2018})
\end{barticle}
\endbibitem

\bibitem[\protect\citeauthoryear{Geoffrion}{2009}]{geoffrion2009lagrangean}
\begin{bchapter}
\bauthor{\bsnm{Geoffrion}, \binits{A.M.}}:
\bctitle{Lagrangean relaxation for integer programming}.
In: \bbtitle{Approaches to Integer Programming},
pp. \bfpage{82}--\blpage{114}.
\bpublisher{Springer},
\blocation{Amsterdam, Netherlands}
(\byear{2009})
\end{bchapter}
\endbibitem

\bibitem[\protect\citeauthoryear{Hiriart-Urruty and Lemar{\'e}chal}{2004}]{hiriart2004fundamentals}
\begin{bbook}
\bauthor{\bsnm{Hiriart-Urruty}, \binits{J.-B.}},
\bauthor{\bsnm{Lemar{\'e}chal}, \binits{C.}}:
\bbtitle{Fundamentals of Convex Analysis}.
\bpublisher{Springer},
\blocation{New York, USA}
(\byear{2004})
\end{bbook}
\endbibitem

\bibitem[\protect\citeauthoryear{Dey et~al.}{2025}]{dey2025geoffrion}
\begin{botherref}
\oauthor{\bsnm{Dey}, \binits{S.S.}},
\oauthor{\bsnm{Meunier}, \binits{F.}},
\oauthor{\bsnm{Ramirez}, \binits{D.M.}}:
Geoffrion's theorem beyond finiteness and rationality.
arXiv preprint arXiv:2510.10966
(2025)
\end{botherref}
\endbibitem

\bibitem[\protect\citeauthoryear{Lemar{\'e}chal and Renaud}{2001}]{lemarechal2001geometric}
\begin{barticle}
\bauthor{\bsnm{Lemar{\'e}chal}, \binits{C.}},
\bauthor{\bsnm{Renaud}, \binits{A.}}:
\batitle{A geometric study of duality gaps, with applications}.
\bjtitle{Mathematical Programming}
\bvolume{90}(\bissue{3}),
\bfpage{399}--\blpage{427}
(\byear{2001})
\end{barticle}
\endbibitem

\bibitem[\protect\citeauthoryear{Dubois-Taine and d’Aspremont}{2025}]{dubois2025frank}
\begin{botherref}
\oauthor{\bsnm{Dubois-Taine}, \binits{B.}},
\oauthor{\bsnm{d’Aspremont}, \binits{A.}}:
Frank-wolfe meets shapley-folkman: a systematic approach for solving nonconvex separable problems with linear constraints: B. dubois-taine, a. d’aspremont.
Mathematical Programming,
1--51
(2025)
\end{botherref}
\endbibitem

\bibitem[\protect\citeauthoryear{Molinaro}{2023}]{molinaro2023strong}
\begin{barticle}
\bauthor{\bsnm{Molinaro}, \binits{M.}}:
\batitle{Strong convexity of feasible sets in off-line and online optimization}.
\bjtitle{Mathematics of Operations Research}
\bvolume{48}(\bissue{2}),
\bfpage{865}--\blpage{884}
(\byear{2023})
\end{barticle}
\endbibitem

\bibitem[\protect\citeauthoryear{Rockafellar}{2015}]{rockafellar2015convex}
\begin{botherref}
\oauthor{\bsnm{Rockafellar}, \binits{R.T.}}:
Convex analysis:(pms-28)
(2015)
\end{botherref}
\endbibitem

\bibitem[\protect\citeauthoryear{Liu and Grimmer}{2023}]{liu2023gauges}
\begin{botherref}
\oauthor{\bsnm{Liu}, \binits{N.}},
\oauthor{\bsnm{Grimmer}, \binits{B.}}:
Gauges and accelerated optimization over smooth and/or strongly convex sets.
arXiv preprint arXiv:2303.05037
(2023)
\end{botherref}
\endbibitem

\end{thebibliography}

\end{document}